\documentclass{article}

\usepackage{cite}
\usepackage[T1]{fontenc}
\usepackage{amsmath,amssymb,amsfonts}
\usepackage{algorithmic,algorithm}
\usepackage{graphicx}
\usepackage{textcomp}
\usepackage{xcolor}
\usepackage{url}
\usepackage{hyperref}
\usepackage{wrapfig}

\renewcommand{\div}{\nabla \cdot}

\newcommand{\OT}{{\Omega_T}}

\newcommand{\GammaT}{{\Gamma_T}}

\newcommand{\ptt}{\partial_{tt}}
\newcommand{\pt}{\partial_{t}}

\newcommand{\pn}{\partial_{n}}

\def\BibTeX{{\rm B\kern-.05em{\sc i\kern-.025em b}\kern-.08em
    T\kern-.1667em\lower.7ex\hbox{E}\kern-.125emX}}
    
\begin{document}

\author{Eric Lindström, \\ 
\footnotesize \textit{Department of Mathematical Sciences}\\
\footnotesize \textit{Chalmers University of Technology}\\
\footnotesize \textit{and University of Gothenburg,}\\
\footnotesize SE-42196 Gothenburg, Sweden, \\
\footnotesize erilinds@chalmers.se
\and
Larisa Beilina, \\ 
\footnotesize \textit{Department of Mathematical Sciences}\\
\footnotesize \textit{Chalmers University of Technology}\\
\footnotesize \textit{and University of Gothenburg,}\\
\footnotesize SE-42196 Gothenburg, Sweden, \\\footnotesize larisa.beilina@chalmers.se}

\title{\Large
Finite element 3D models of melanoma growth and time-dependent backscattered data for dielectric properties of melanoma at 6 GHz
}
\date{\small \today}

\maketitle

\vspace{-.7cm}
\begin{abstract}

\footnotesize

Finite element meshes for 3D models simulating realistic malignant
melanoma (MM) growth, incorporating accurate dielectric properties of
the skin, have been developed. Numerical simulations illustrate how 3D
finite element meshes can be utilized to generate backscattered data,
enabling the evaluation of reconstruction algorithms designed to
determine the dielectric properties of the proposed 3D model.
 
\end{abstract}



\graphicspath{
  {FIGURES/}
  {/chalmers/users/larisa/VR_Sweden/VR_2024/PLAN1/FIGURES/}
 {/chalmers/groups/larisa_3/Acoustic/ISIC_sampled100/RecISIC1transfreq40/}
 {pics/}}


\section{Introduction}

This study presents finite element meshes for 3D models that simulate
the growth of realistic malignant melanoma (MM). These models are
integrated with accurate representations of skin properties
corresponding to the regions where MM develops.  To generate finite
element meshes we combine geometrical models of melanoma growth
presented recently in \cite{IEEE2024} with skin models studied in
\cite{ICEAA2024_BEN}.

One of the key applications of the developed models is microwave
medical imaging \cite{ref1}   for detecting malignant melanoma, the most lethal
form of skin cancer, despite accounting for only 1\% of all skin cancer
cases \cite{IEEE2024,ref1}.
The prognosis for MM is closely linked to the depth of primary
tumor invasion in the skin  \cite{skindepth,ICEAA2024_BEN,IEEE2024}. Microwave medical imaging stands out as a
non-invasive imaging technique, offering a significant advantage over
methods such as X-ray imaging   \cite{noninvasive,Pastorino}. This makes it a compelling addition to
existing medical imaging modalities like ultrasound  \cite{gonch1,gonch2}, X-ray, and
MRI \cite{tomography}. Notably, X-ray and MRI are not utilized for diagnosing primary
skin cancers. Consequently, microwave imaging emerges as an attractive
non-invasive technology that can be applied directly to the
skin \cite{IEEE2024}. Developing this technology holds great promise as a valuable
complement to current diagnostic approaches for MM.

Previous studies \cite{IEEE2024,ieee1,wisconsin}
have highlighted variations in permittivity contrasts
between malignant and normal tissues, demonstrating that malignant
tumors exhibit higher relative permittivity values compared to normal
tissues at frequencies below 10 GHz. However, accurately
estimating the relative permittivity of the skin's internal structure
with malignant melanoma (MM) remains a significant challenge. This
estimation relies on analyzing backscattered electromagnetic waves
collected by detectors positioned near or on the skin affected by
MM. In this study, we have modeled MM using realistic dielectric
properties measured at 6 GHz, based on data from \cite{IEEE2024}.

Potential application of the finite element meshes presented in this
 work is qualitative and quantitative malign melanoma detection using
 microwaves.  One possibility to do this is to use an adaptive domain
 decomposition finite element/finite difference method (ADDFE/FDM)
 developed recently in \cite{BL1}.  This method was applied in
 \cite{BL1,BL2,LB4} for reconstruction of dielectric permittivity and
 conductivity functions at 6 GHz for anatomically realistic breast
 phantom from database of \cite{wisconsin}, as well as in
 \cite{ICEAA2025_KLB} for reconstructing the dielectric properties of
 MM for another halph-spheric 3D melanoma model.

Additional challenge is apply real measured backscattered microwave
data at different frequencies in order to use generated finite element meshes in the ADDFE/FDM
method for qualitative and quantitative MM detection. This challenge can be considered as a future work which should be performed in collaboration with electrical engineers.

 The problem
 of determination of the spatially distributed relative dielectric
 permittivity and conductivity functions using backscattered data of
 the electric field collected at the boundary of the investigated
 domain  in scientific community is called Coefficient
 Inverse Problems (CIPs).
  We refer to 
 \cite{BakKok},\cite{BK}, \cite{Ghavent}, \cite{GG},  \cite{T},  \cite{itojin},  \cite{KSS} and
 references therein  for mathematical methods related to the solution of CIPs  using  minimization of the functional,
 and to \cite{BK,BuKr,convex1,TBKF2,KBKSNF}  for 
 globally convergent methods for solution of CIPs.

In the numerical studies of the current work  we demonstrate how
the developed 3D finite element mesh with MM can be applied in the ADDFE/FDM
method of \cite{BL1} to generate backscattered data. This data serves as a testing
platform for reconstruction algorithms aimed at determining the
dielectric properties of the proposed 3D models.

\begin{wrapfigure}{r}{0.4\textwidth}
\begin{center}
    \includegraphics[trim={1cm 2cm 0cm 1.5cm},clip,width=1.05\linewidth]{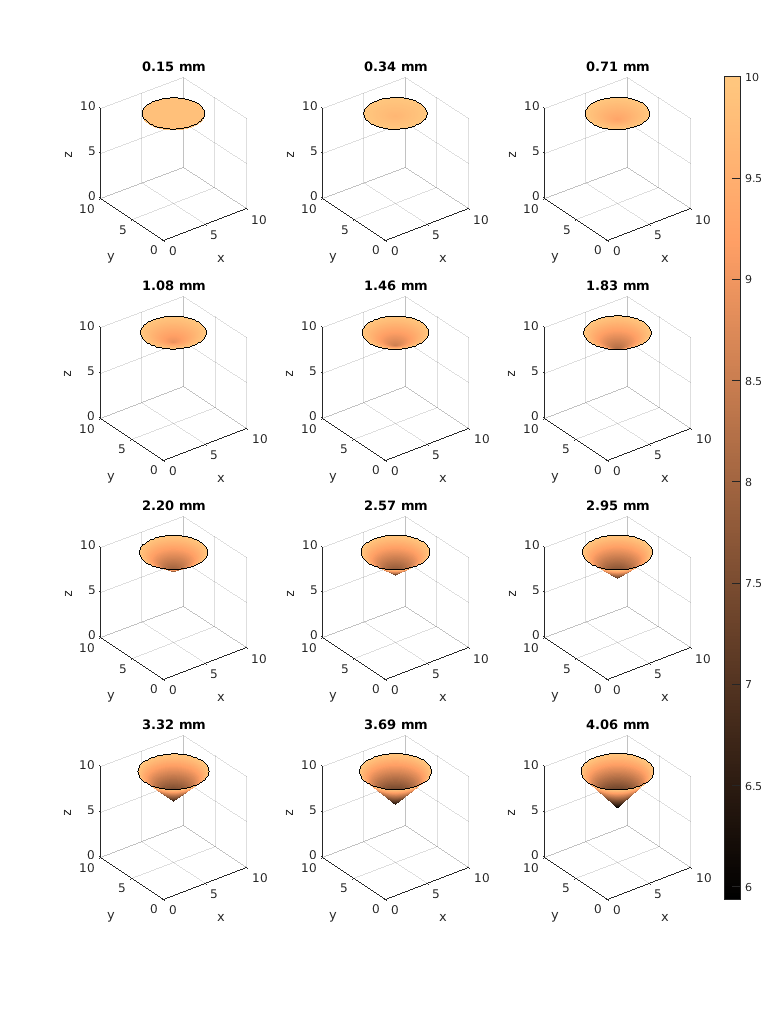}
    \caption{\textit{Melanoma model shapes for different depths.}}
    \label{fig:GeomMM}
\end{center}
\end{wrapfigure}

\section{Melanoma growth models}

\label{sec:MM}

In this study we are taking geometrical models of MM growth developed
 in \cite{IEEE2024} and combine them with accurate representation of
 skin properties as was proposed in \cite{ICEAA2024_BEN}.  The
 prognosis of MM is related to the diameter and depth of MM invasion \cite{28,29}
 and we modelled geometrical properties and shapes of MM accordingly
 to Figure \ref{fig:GeomMM} and Tables 1, 2.  Next, we combined models
 shown on Figure \ref{fig:GeomMM} with skin properties described in
 the Table \ref{tab:table1}.  We note that this table highlights the
 dielectric properties of skin with MM at 6 GHz. However, the depth of
 all tissue types are consistent across all frequencies.

Let  the computational domain  $\Omega\subset \mathbb{R}^3$ be a convex domain with a smooth boundary $\Gamma$.
To generate finite element meshes we choose  $\Omega$ of the size $10 \times 10 \times 10 $ mm   such that
the dimensionless domain is
 \begin{equation*}
    \Omega = \left\{ x= (x_1,x_2, x_3) \in
           (0,  10)
    \times (0,  10)
    \times (0, 10)
    \right\}.
 \end{equation*}
  Next, we decompose $\Omega$ into tetrahedral elements with the mesh size $h=0.5$
 such that every tetrahedral element has 
type of material corresponding to the tissue type of the Tables \ref{fig:GeomMM},
\ref{tab:table1},
 and assign then values for the relative dielectric permittivity $\varepsilon_r$ and effective conductivity
   $\sigma$ accordingly to the Tables \ref{fig:GeomMM},
\ref{tab:table1}.

\begin{figure*}
    \centering
    \begin{tabular}{c c c c}
         Month 0 & Month 2 & Month 4 & Month 6 \\
         \includegraphics[trim={1cm 0cm 3cm 2.5cm},clip,width=0.21\linewidth]{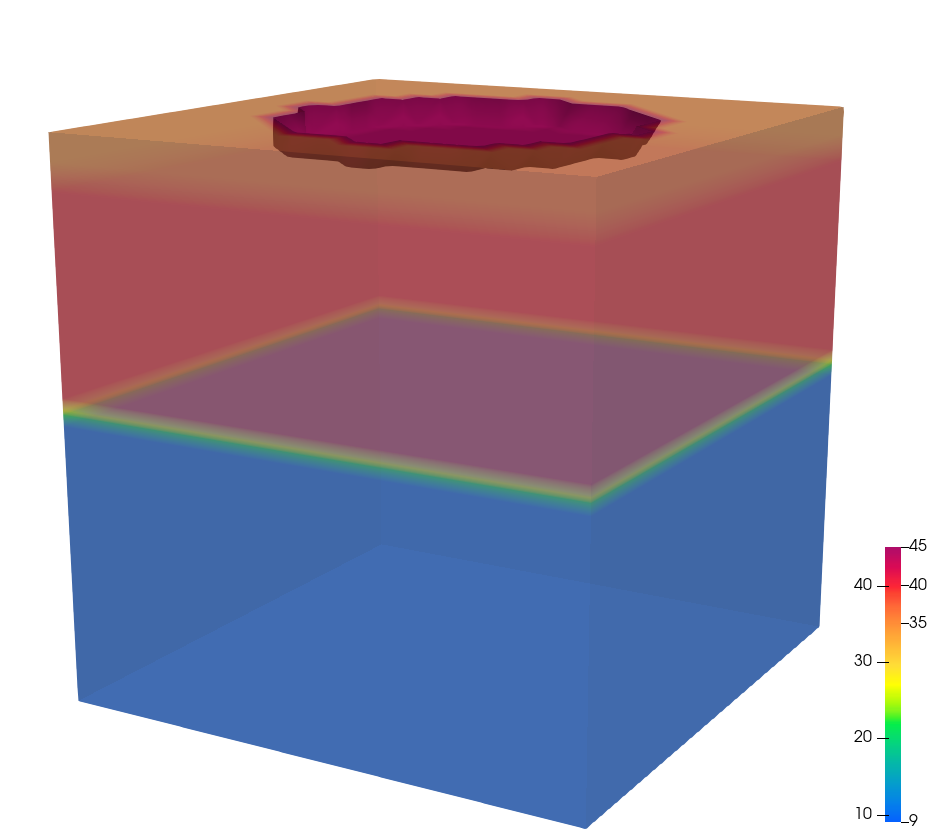} & \includegraphics[trim={1cm 0cm 3cm 2.5cm},clip,width=0.21\linewidth]{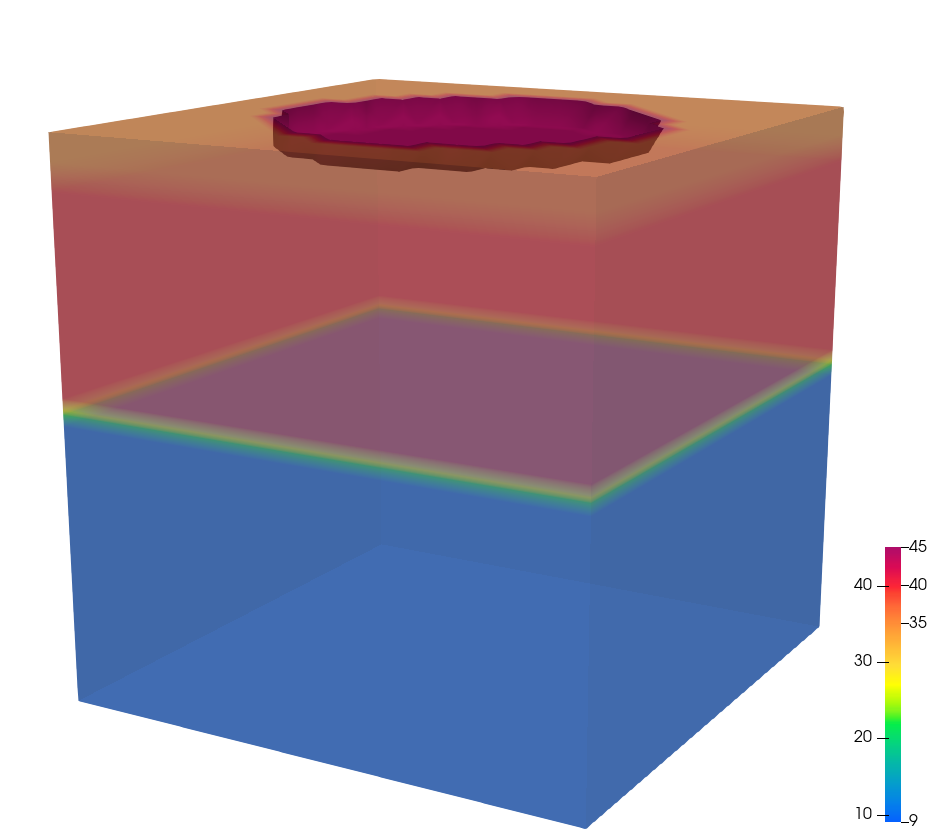} & \includegraphics[trim={1cm 0cm 3cm 2.5cm},clip,width=0.21\linewidth]{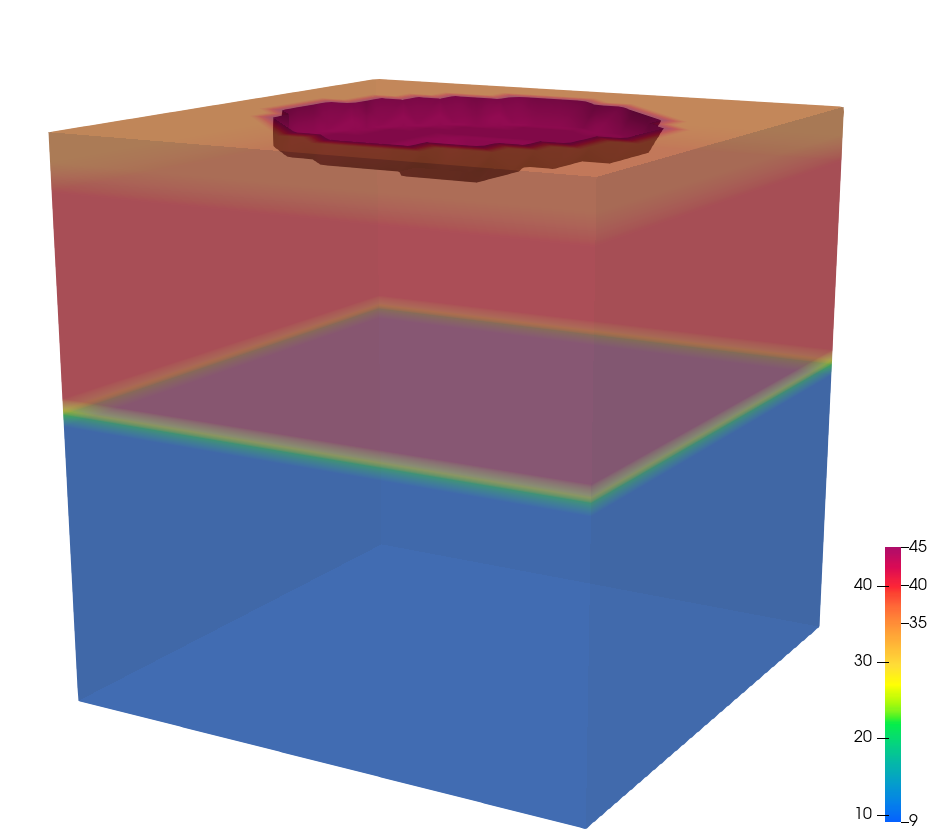} & \includegraphics[trim={1cm 0cm 3cm 2.5cm},clip,width=0.21\linewidth]{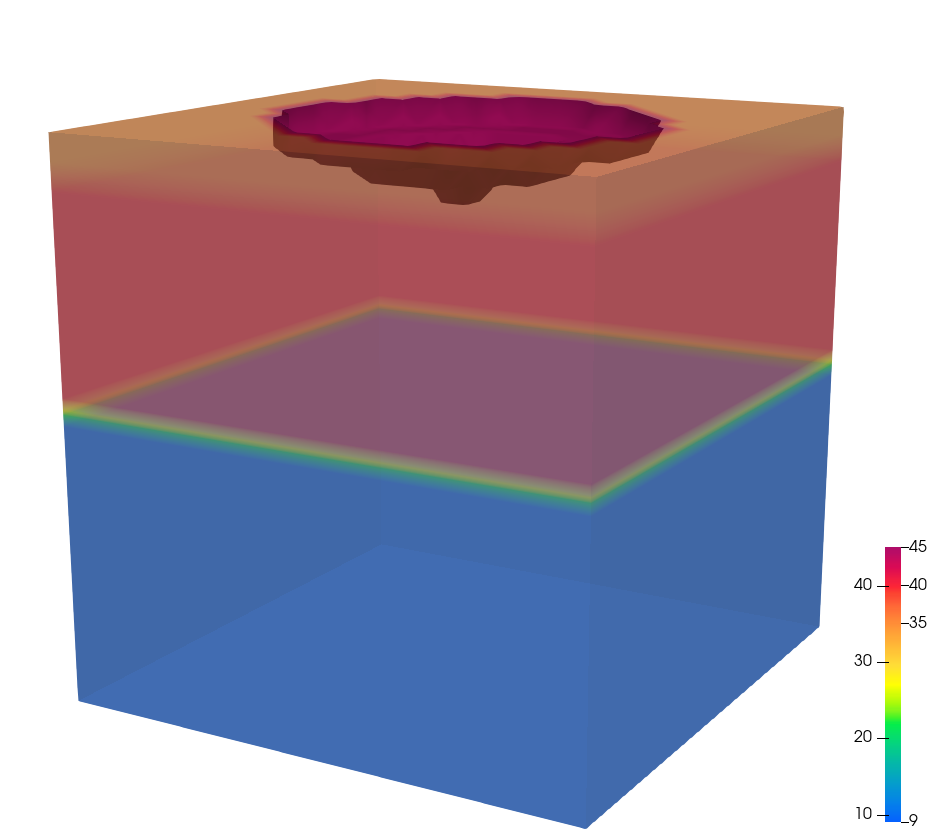} \\
         Month 8 & Month 10 & Month 12 & Month 14 \\
         \includegraphics[trim={1cm 0cm 3cm 2.5cm},clip,width=0.21\linewidth]{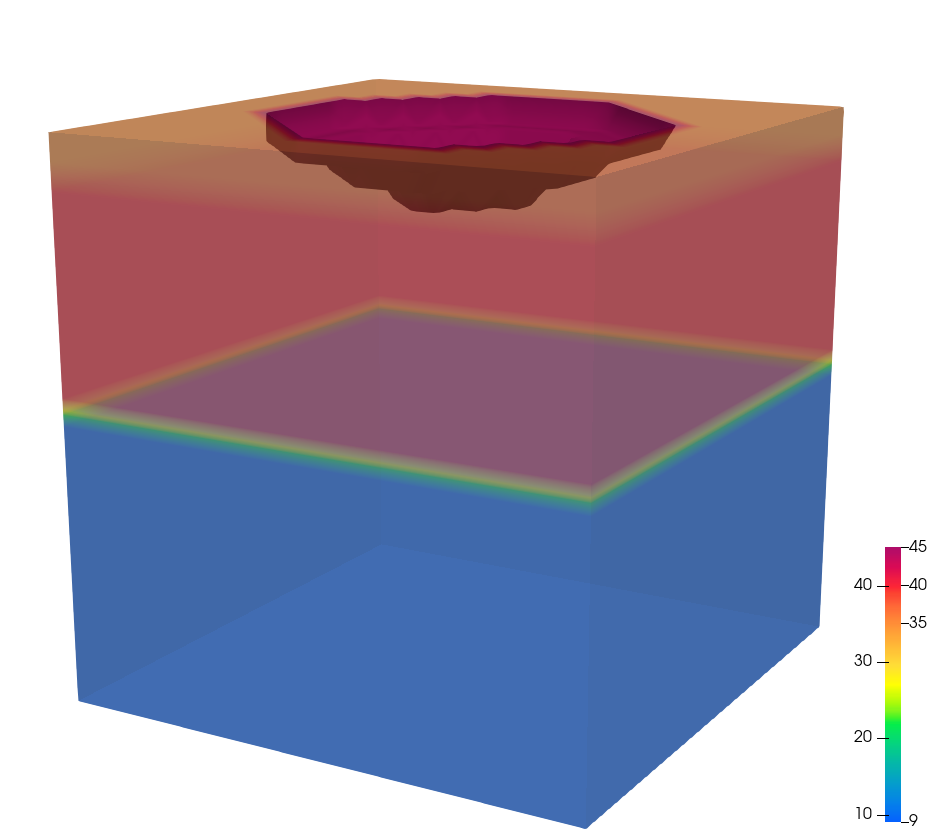} & \includegraphics[trim={1cm 0cm 3cm 2.5cm},clip,width=0.21\linewidth]{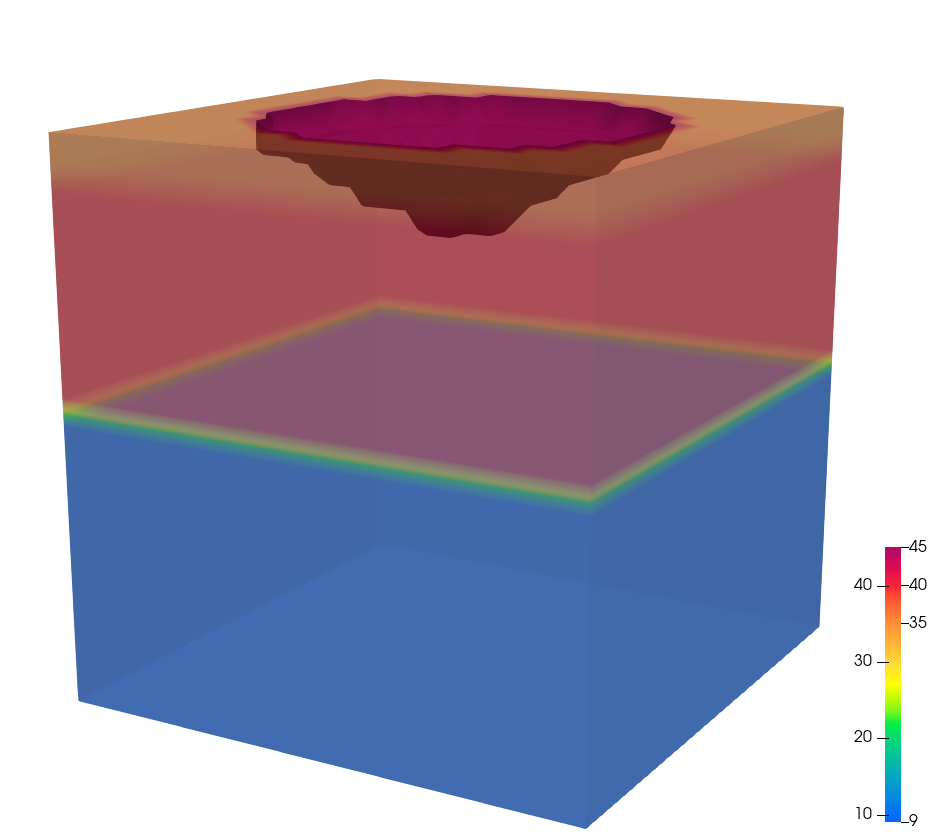} & \includegraphics[trim={1cm 0cm 3cm 2.5cm},clip,width=0.21\linewidth]{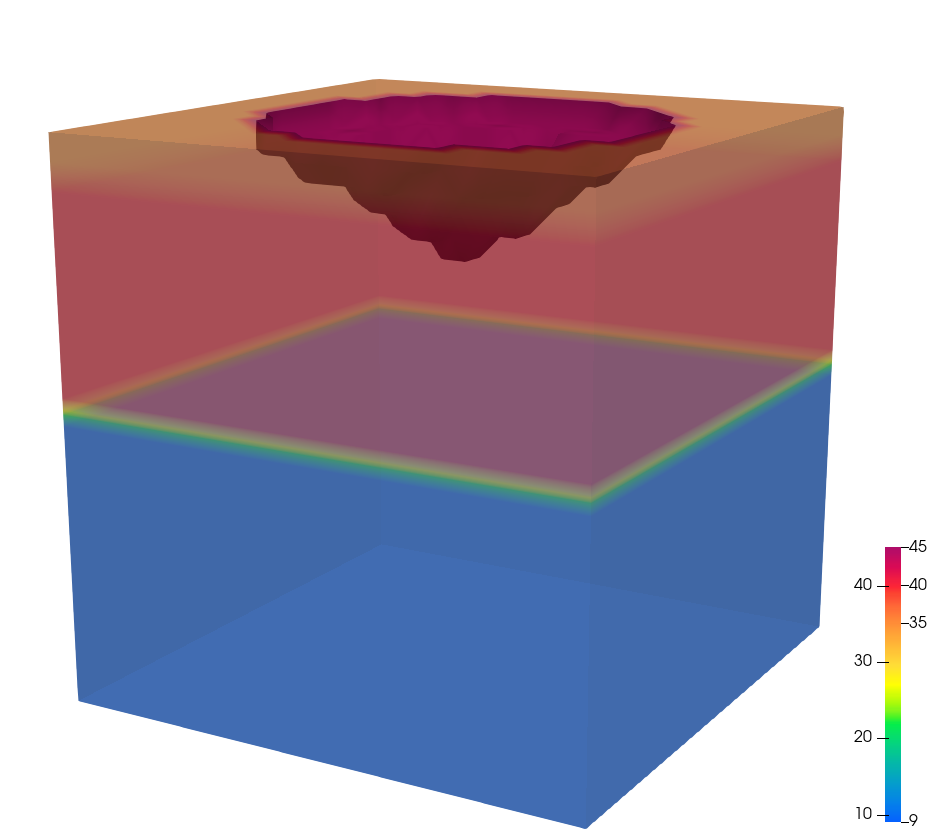} & \includegraphics[trim={1cm 0cm 3cm 2.5cm},clip,width=0.21\linewidth]{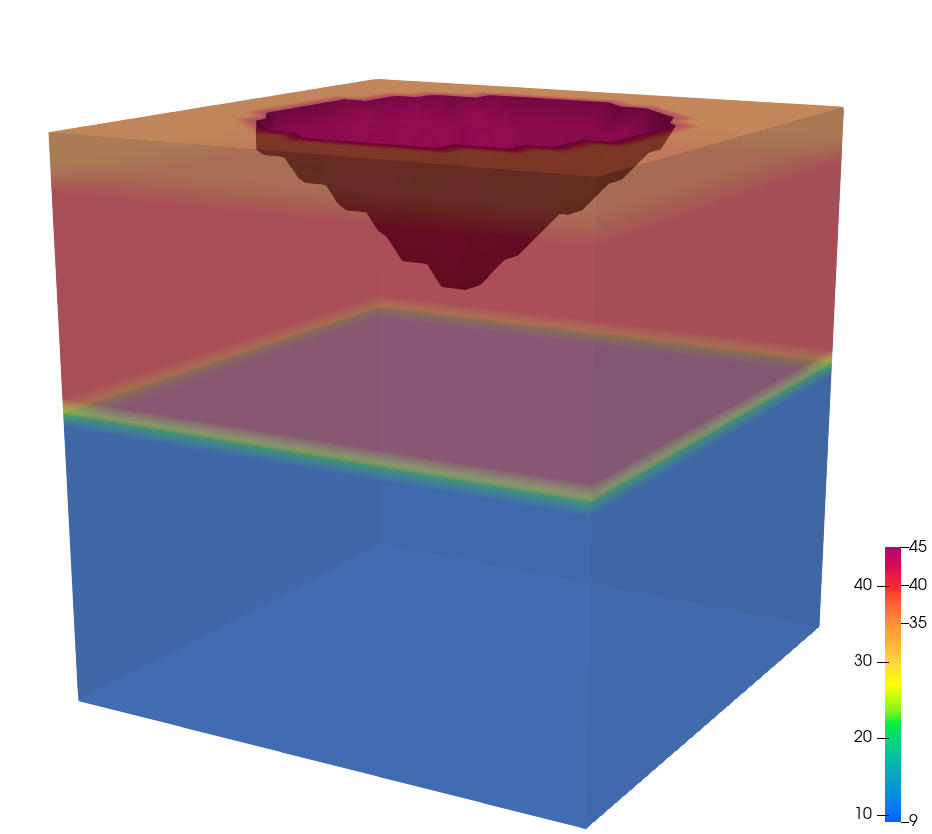} \\
         Month 16 & Month 18 & Month 20 & Month 22 \\
         \includegraphics[trim={1cm 0cm 3cm 2.5cm},clip,width=0.21\linewidth]{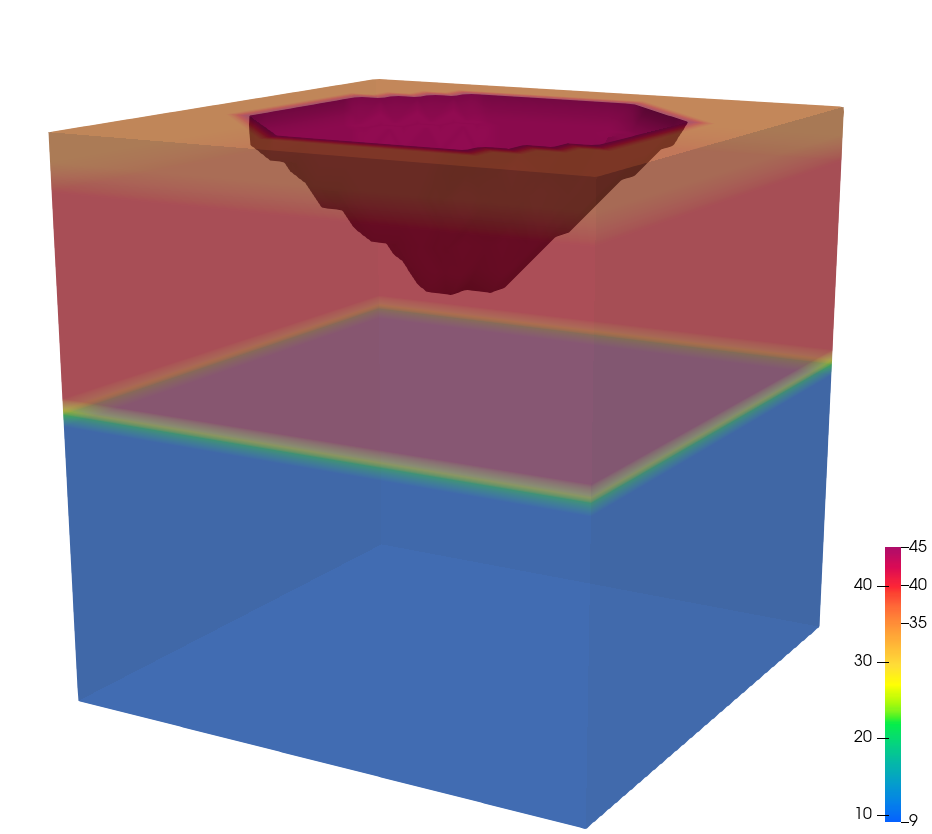} & \includegraphics[trim={1cm 0cm 3cm 2.5cm},clip,width=0.21\linewidth]{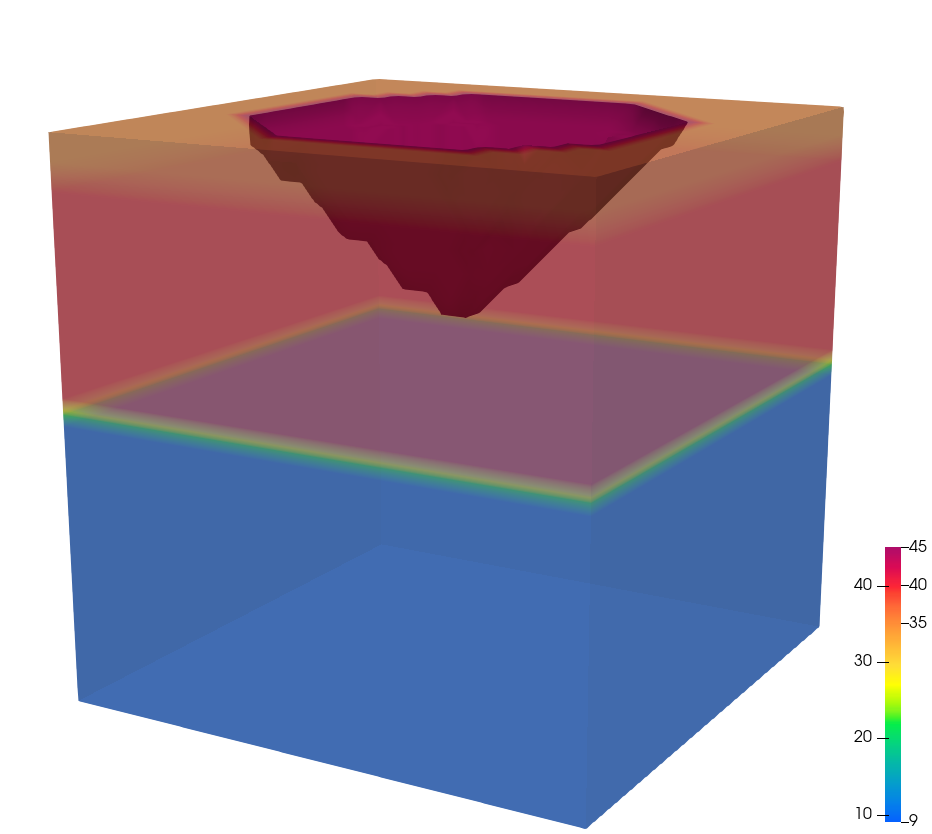} & \includegraphics[trim={1cm 0cm 3cm 2.5cm},clip,width=0.21\linewidth]{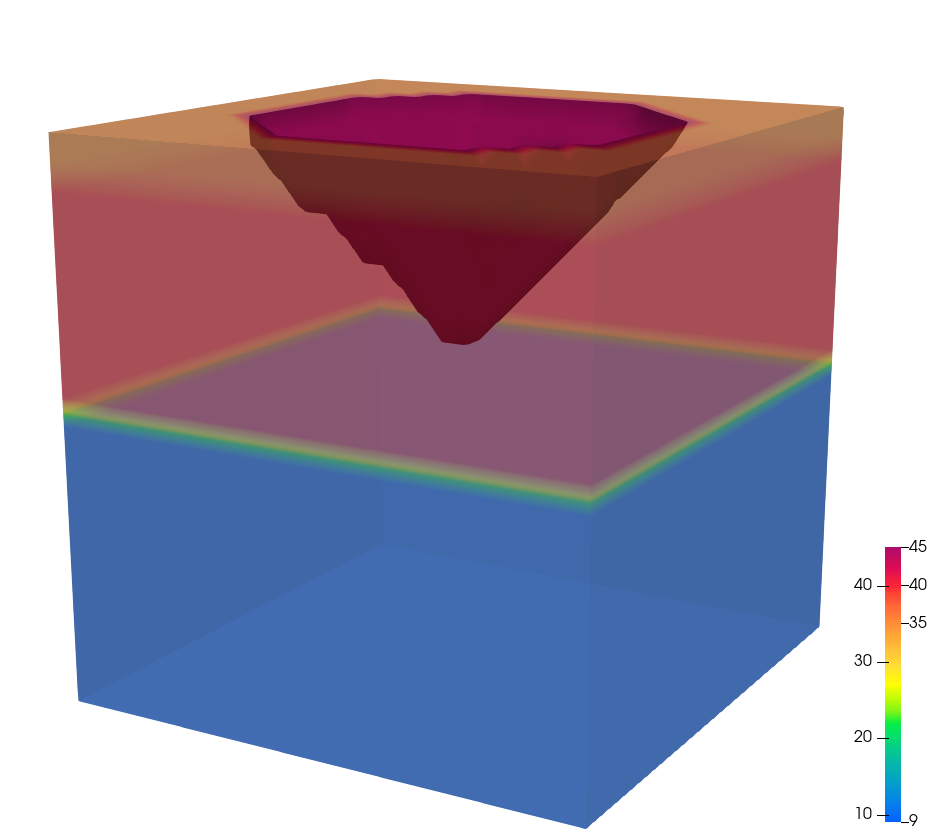} & \includegraphics[trim={1cm 0cm 3cm 2.5cm},clip,width=0.21\linewidth]{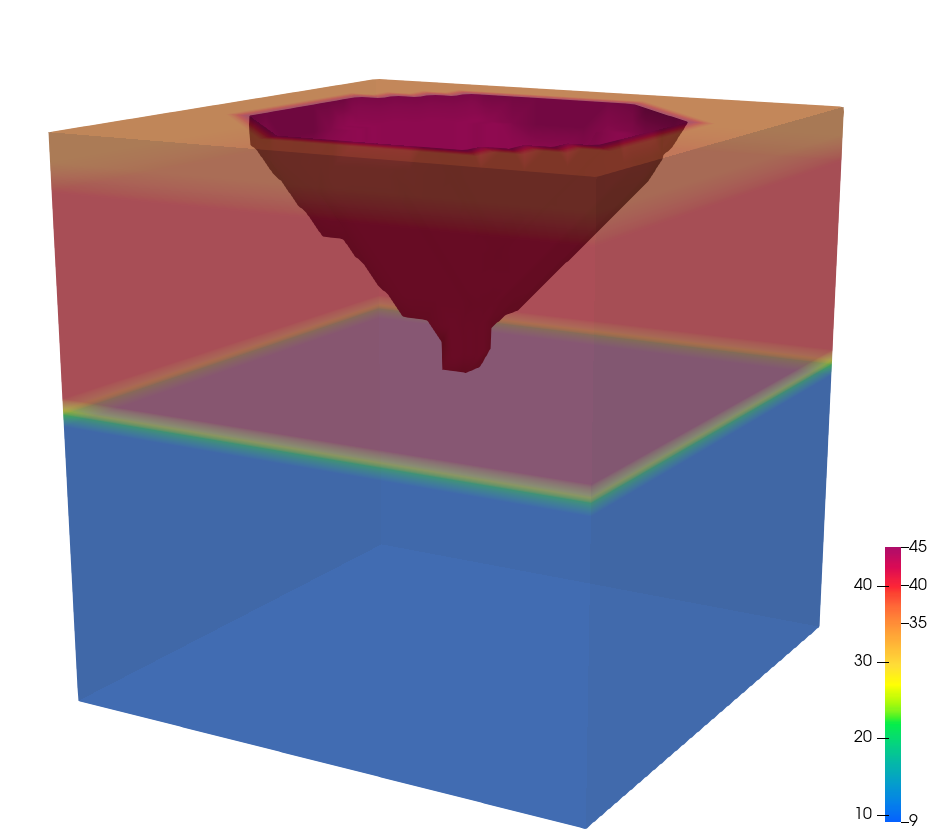} \\
    \end{tabular}
    \caption{\textit{Models used for calculations visualized based on geometric properties shown in \autoref{geotab}.}}
    \label{fig:model3D}
\end{figure*}

 Figure
  \ref{fig:model3D} presents finite element models of MM  corresponding to  the geometrical models presented in Figure  \ref{fig:GeomMM} and Table  \ref{geotab} combined
  with realistic values of skin shown in Table  \ref{tab:table1}.

Using Figure  \ref{fig:model3D}  and Tables \ref{geotab}, \ref{tab:table1}  we observe that MM models for
 months $0,2,4,6,8,10$  correspond to MM "in-situ" \cite{28,29} and these MM are not invasive since the transition from non-invasive to invasive MM starts from the diameter  of MM $d =6$ mm \cite{IEEE2024}. 
Next, MM models for months $12,14,16,18$ are located in epidermis and dermis and they are invasive \cite{30}.
 Finally, MM models for month $20,22$ are propagated from dermis to fat and can reach the lymphatic system where  MM begins to metastasize \cite{31}.

\begin{figure*}
    \centering
    \begin{tabular}{c c c c}
         Month 0 & Month 2 & Month 4 & Month 6 \\
         \includegraphics[trim={8cm 5cm 8cm 5cm},clip,width=0.19\linewidth]{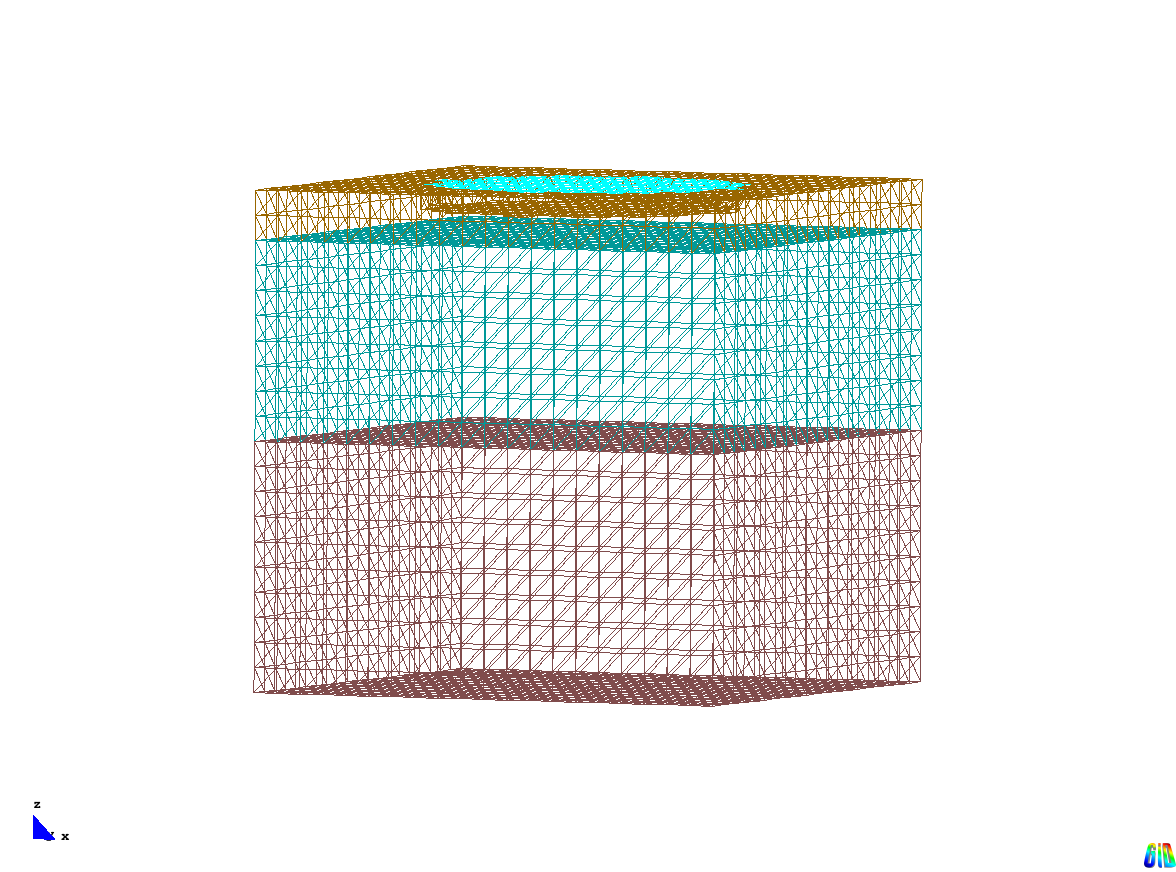} & \includegraphics[trim={8cm 5cm 8cm 5cm},clip,width=0.19\linewidth]{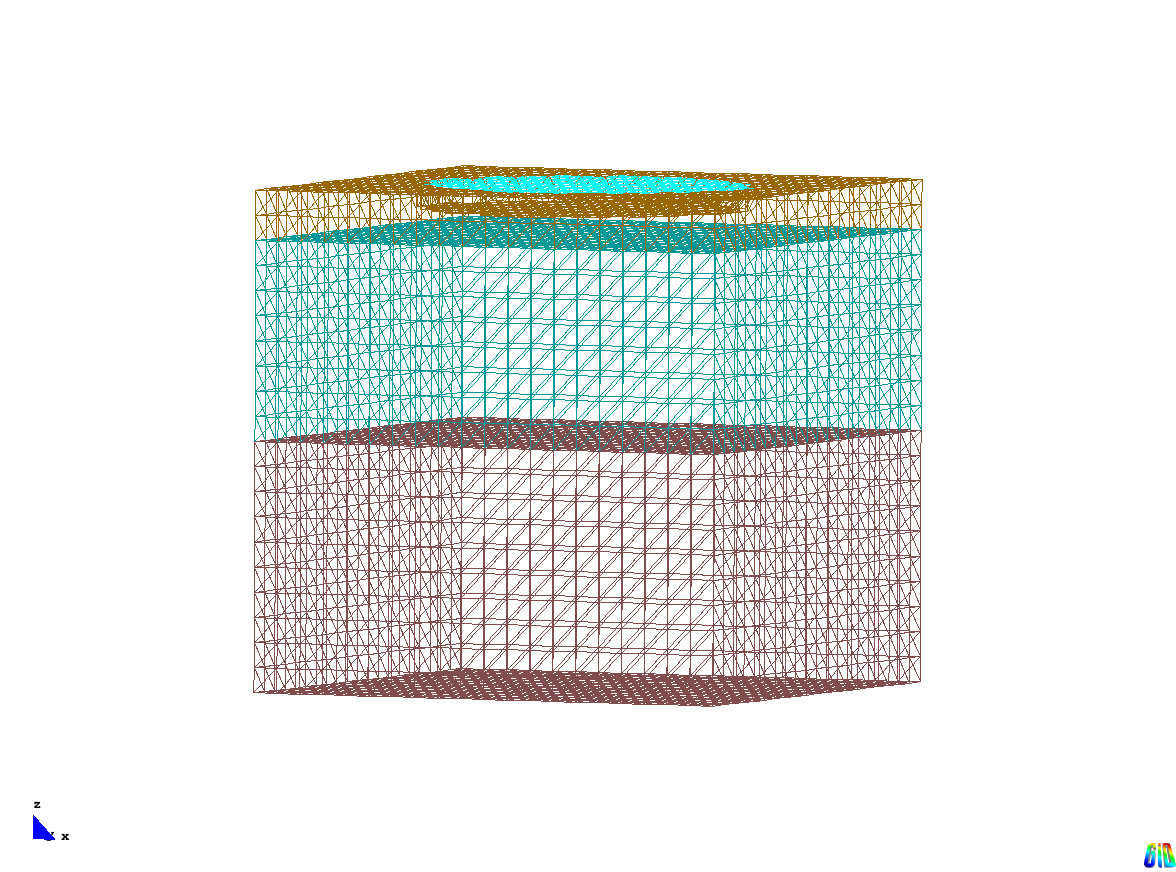} & \includegraphics[trim={8cm 5cm 8cm 5cm},clip,width=0.19\linewidth]{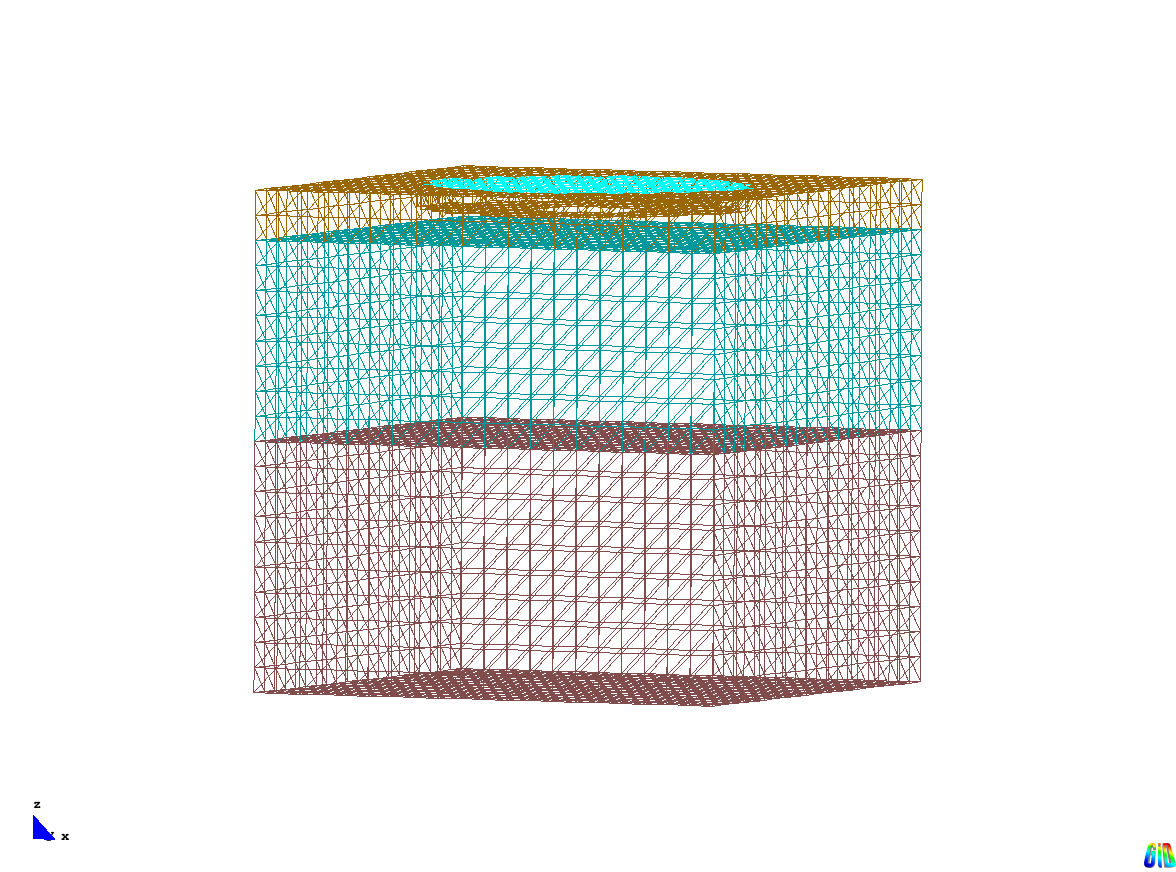} & \includegraphics[trim={8cm 5cm 8cm 5cm},clip,width=0.19\linewidth]{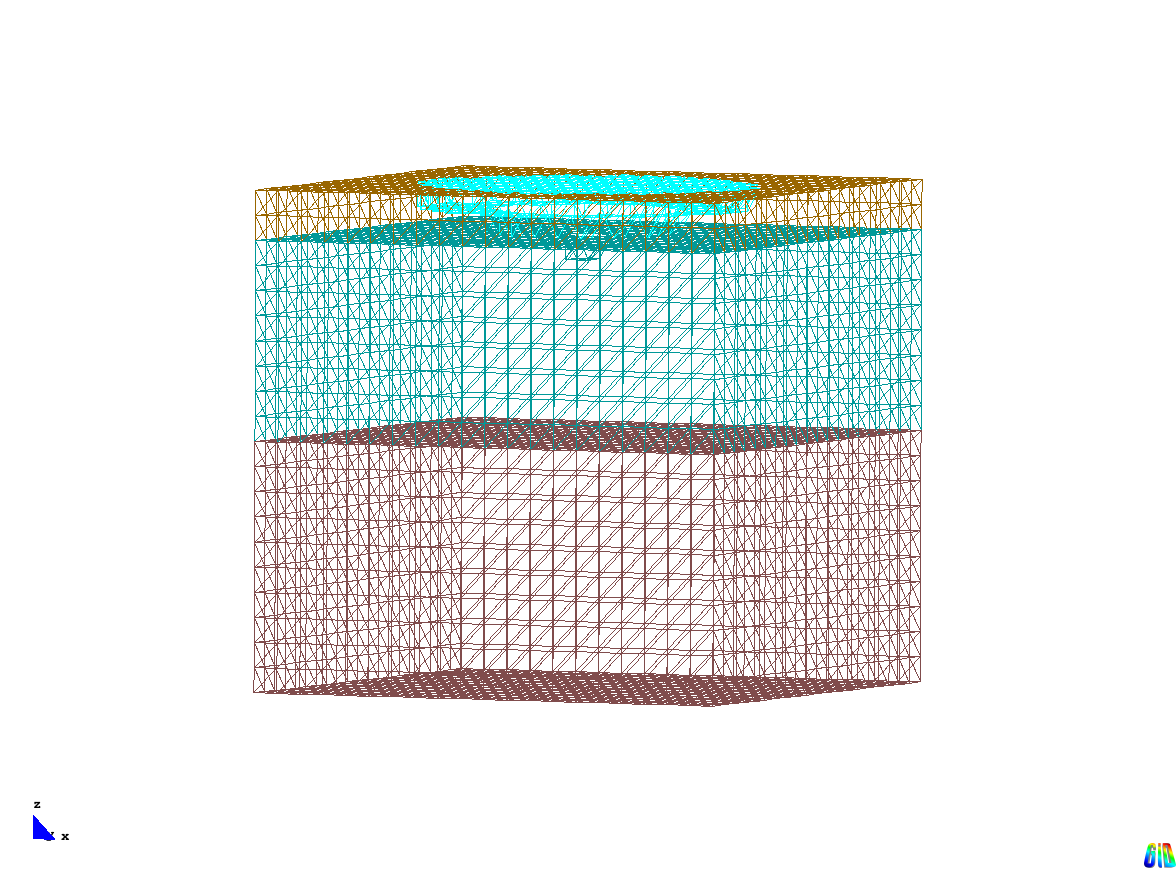} \\
         Month 8 & Month 10 & Month 12 & Month 14 \\
         \includegraphics[trim={8cm 5cm 8cm 5cm},clip,width=0.19\linewidth]{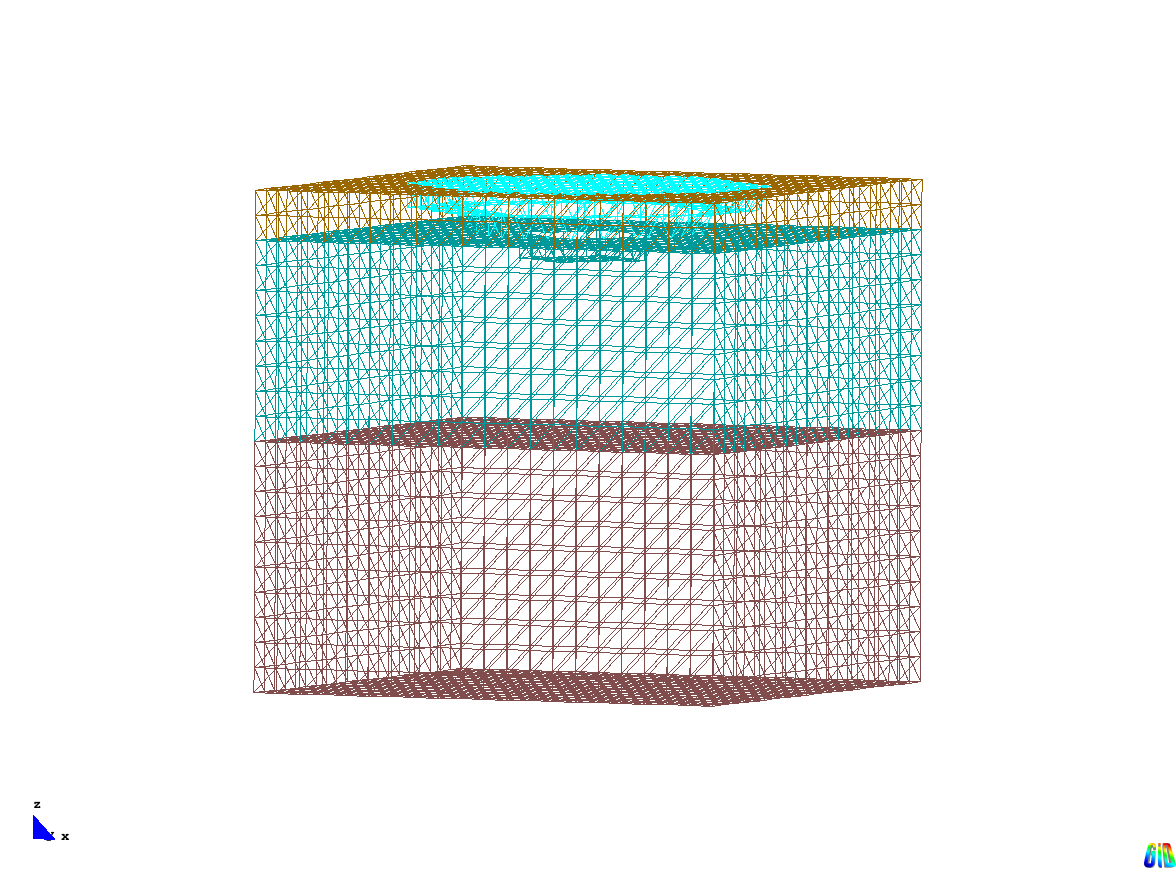} & \includegraphics[trim={8cm 5cm 8cm 5cm},clip,width=0.19\linewidth]{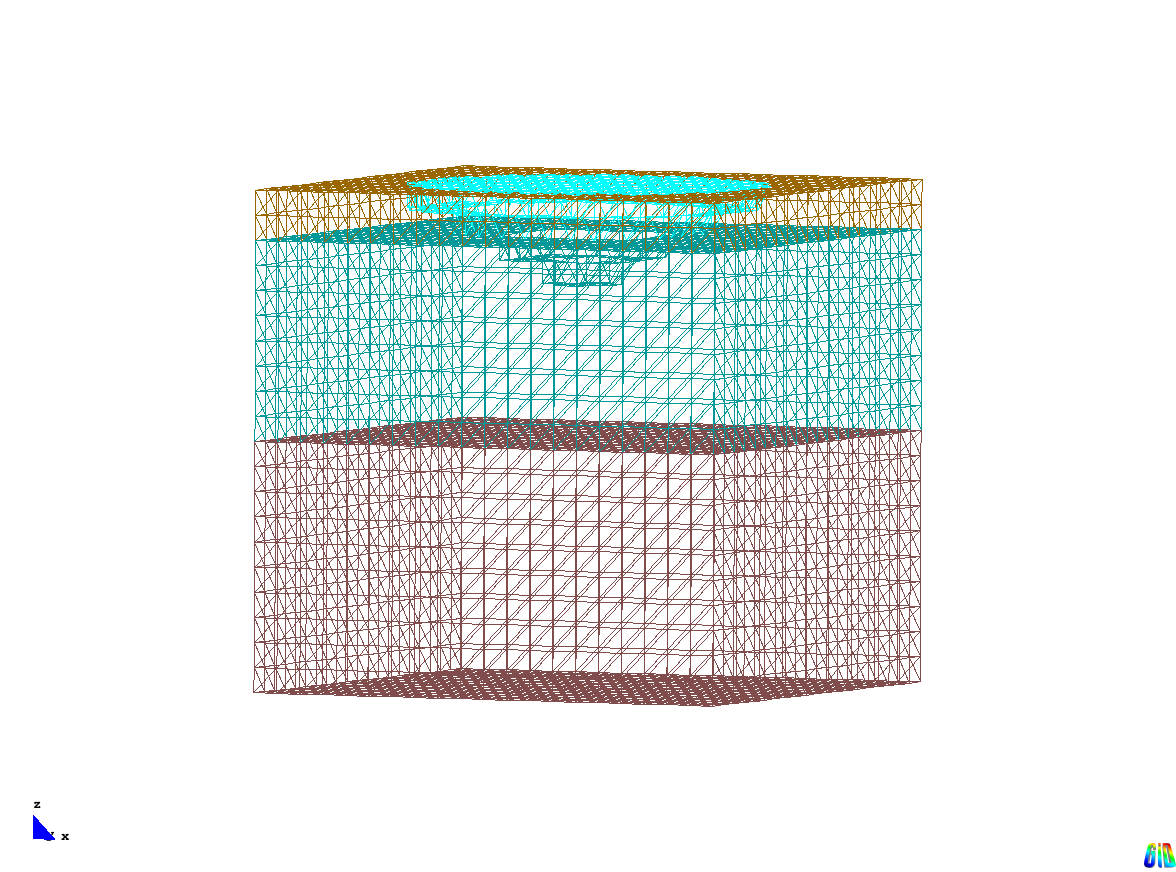} & \includegraphics[trim={8cm 5cm 8cm 5cm},clip,width=0.19\linewidth]{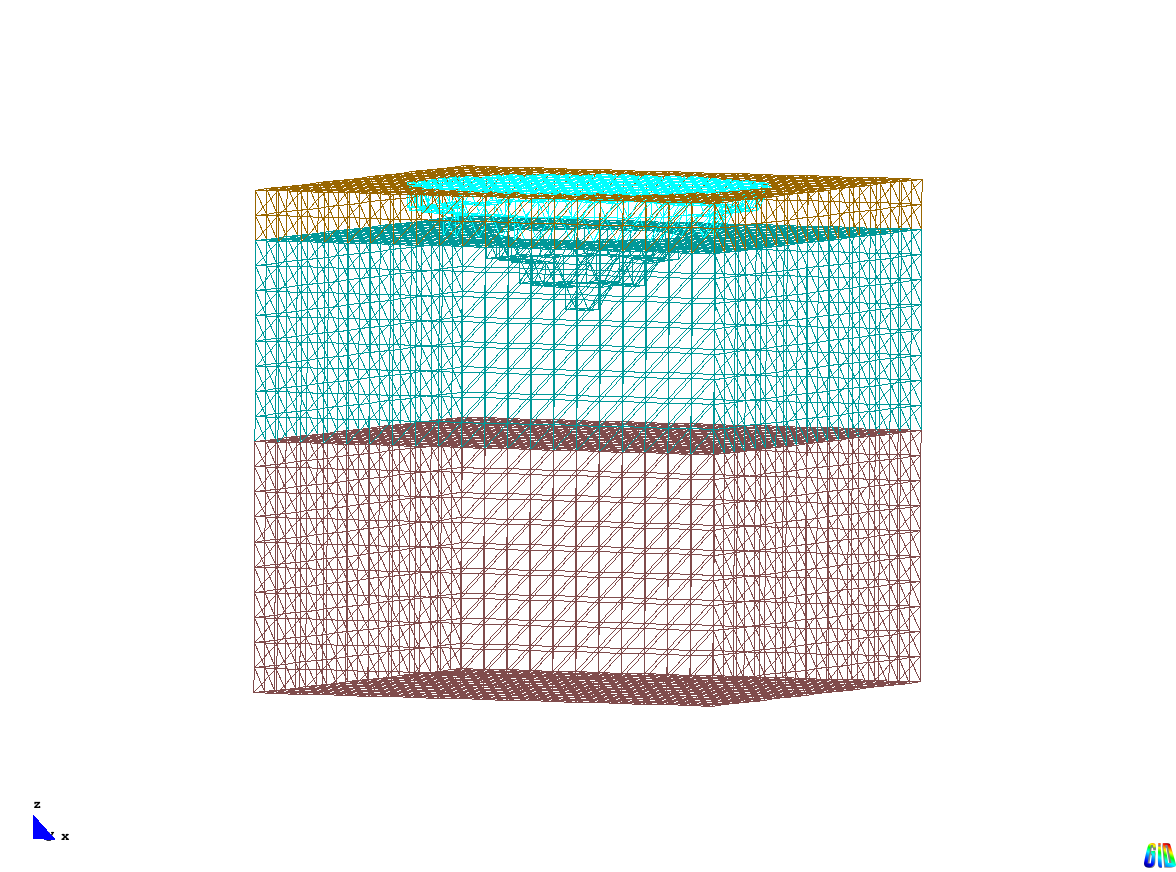} & \includegraphics[trim={8cm 5cm 8cm 5cm},clip,width=0.19\linewidth]{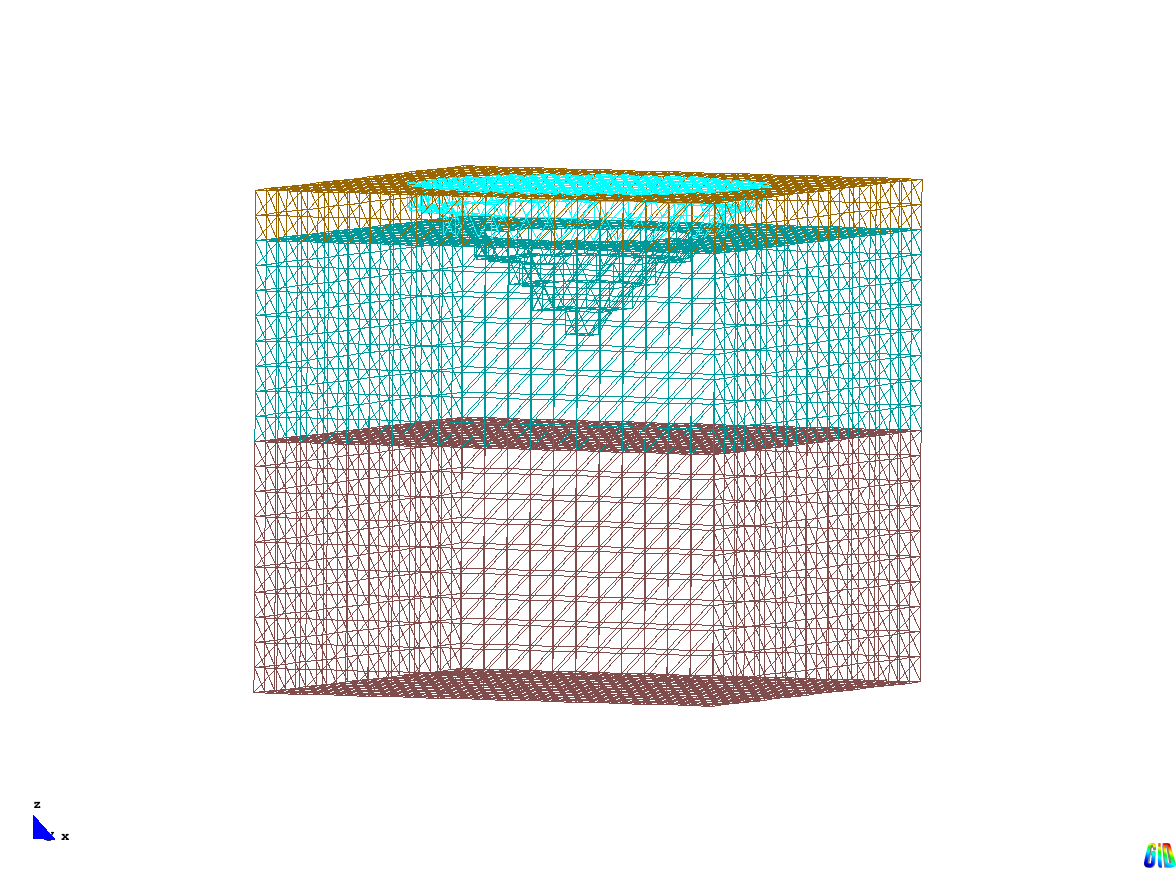} \\
         Month 16 & Month 18 & Month 20 & Month 22 \\
         \includegraphics[trim={8cm 5cm 8cm 5cm},clip,width=0.19\linewidth]{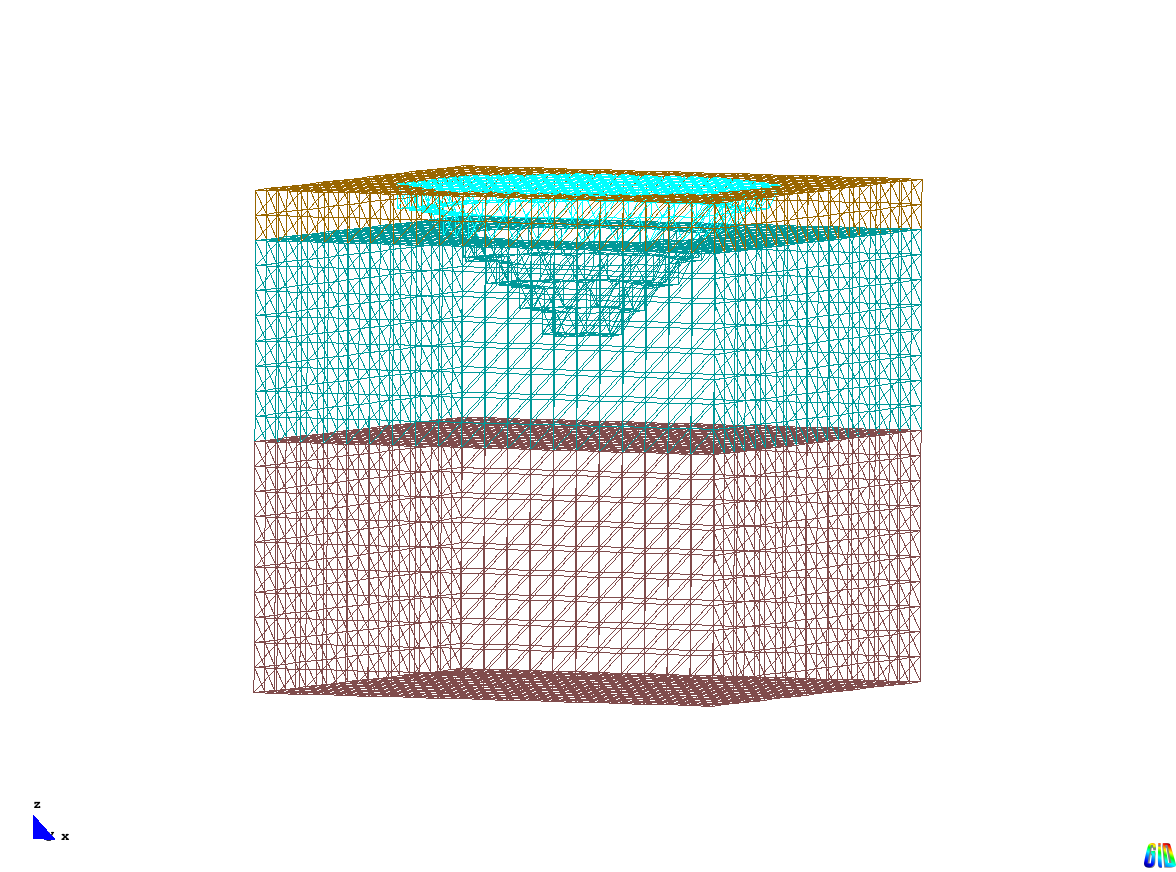} & \includegraphics[trim={8cm 5cm 8cm 5cm},clip,width=0.19\linewidth]{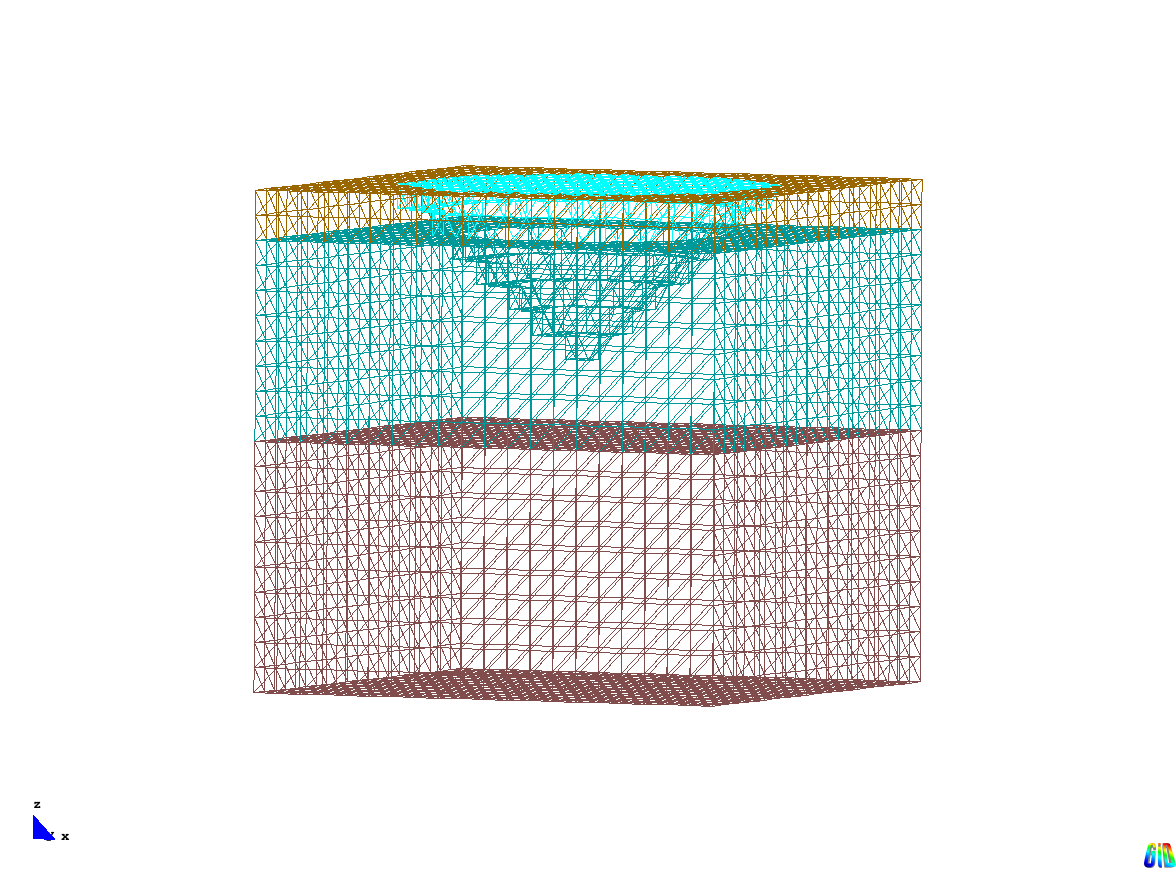} & \includegraphics[trim={8cm 5cm 8cm 5cm},clip,width=0.19\linewidth]{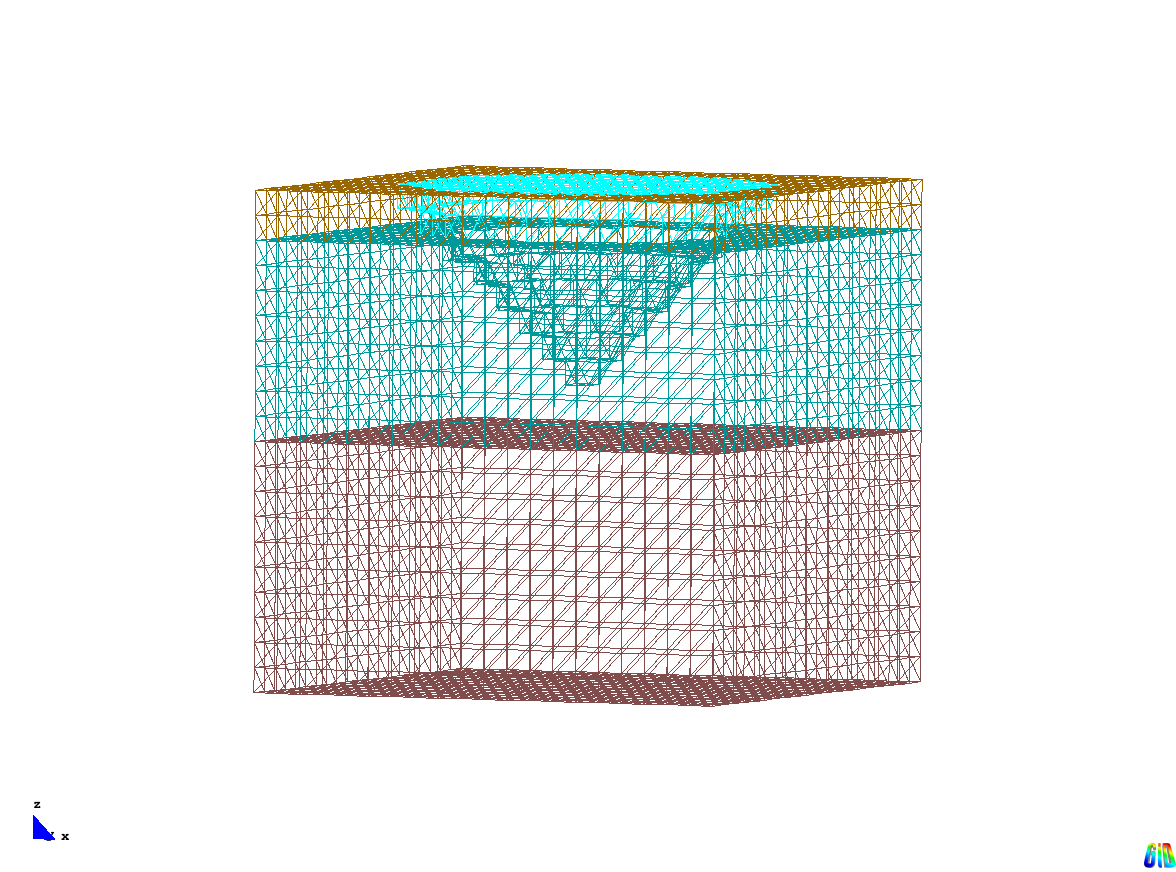} & \includegraphics[trim={8cm 5cm 8cm 5cm},clip,width=0.19\linewidth]{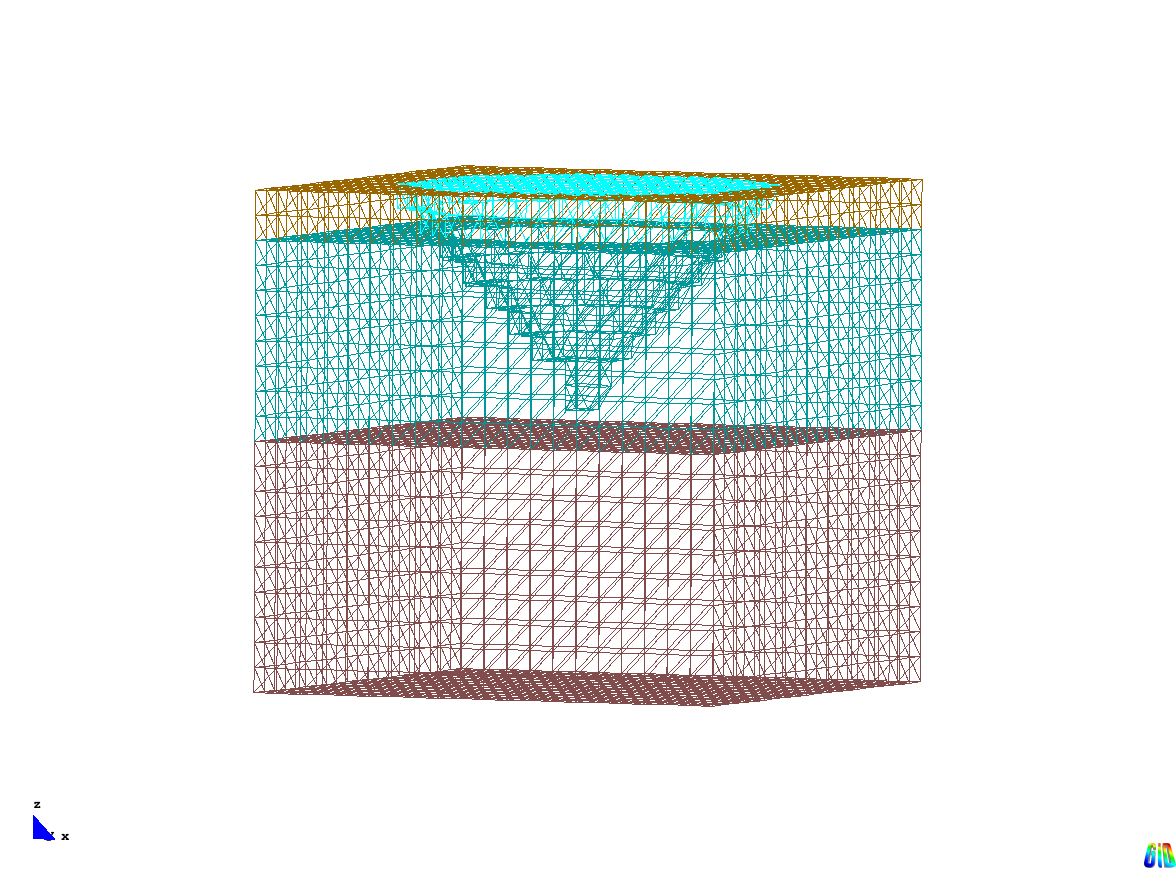} \\
    \end{tabular}
    \caption{\textit{Finite element meshes for models presented on Figure  \ref{fig:model3D}.}}
    \label{fig:meshes3D}
\end{figure*}

\begin{table}[h!]
    \centering
    \vspace{-.5cm}
    \caption{\textit{Tumor model geometrical properties}} 
    \begin{tabular}{|c c c c|}
    \hline 
        Month & Shape & Diameter (mm) & Depth (mm) \\
    \hline
    \hline
        0 & Cylinder & 6 & 0.15 \\
        2 & Cone & 6.09 & 0.34 \\
        4 & Cone & 6.19 & 0.71 \\
        6 & Cone & 6.28 & 1.08 \\
        8 & Cone & 6.36 & 1.46 \\
        10 & Cone & 6.45 & 1.83 \\
        12 & Cone & 6.54 & 2.20 \\
        14 & Cone & 6. 62 & 2.57 \\
        16 & Cone & 6.71 & 2.95 \\
        18 & Cone & 6.79 & 3.32 \\
        20 & Cone & 6.87 & 3.69 \\
        22 & Cone & 6.96 & 4.06 \\
    \hline
    \end{tabular} 
    \label{geotab}
\end{table}

\section{The Mathematical Model}

\label{sec:model}

The mathematical model which we take for simulations in this work is stabilized time-dependent Maxwell's equations for electric field
 for isotropic and linear materials
which approximates original time-dependent Maxwell's equations. This model  was
theoretically   studied in several recent works  - see \cite{BR1,BL1,LB3}.
 Let the electric field $E\left( x,t\right)
=\left( E_{1},E_{2},E_{3}\right) \left( x,t\right), x \in \Omega$ is changing in the time interval $t \in (0, T)$ under the 
assumption  that the dimensionless relative magnetic permeability of 
the medium is $\mu_r \equiv  1$. Taking $\varepsilon = \varepsilon_r$
the stabilized model problem in  $\Omega_T := \Omega \times (0,T)$
is:
\begin{align}
    &\varepsilon \ptt E + \sigma \pt E - \Delta E - \nabla\div (\varepsilon - 1) E = 0  \hspace{-.2cm}&&\text{in } \OT , \label{forweq1}\\
    &E(x, 0) = f_0(x), ~ \pt E(x, 0) = f_1(x) &&\text{in } \Omega  \label{forwinitcond}.
\end{align}

To be able solve the system \eqref{forweq1}--\eqref{forwinitcond}
numerically one need supply it with appropriate boundary conditions.
Usually in real-life scenarios we know values of $\varepsilon, \sigma$
in some part of the domain and these values are not known in the
another part of the domain. In such cases it is convenient to use
domain decomposition of the whole computational domain as it is done
in works \cite{BL1,BL2}.  To be able use the domain decomposition
finite element/finite difference method (DDFE/FDM) of \cite{BL1, BL2}
for the model equation \ref{forweq1} -- \ref{forwinitcond} with finite
element geometries generated in section \ref{sec:MM} we construct the
new hybrid finite element/finite difference computational domain
$\Omega_{\rm hyb}$ such that $\Omega_{\rm hyb} : = \Omega_{\rm FEM}
\cup \Omega_{\rm FDM}$ with $\Omega_{\rm FEM} = \Omega$ and
$\Omega_{\rm FEM} \subset \Omega_{\rm hyb}$. We apply the finite
element method in $\Omega_{\rm FEM}$ and the finite difference method
in $\Omega_{\rm FDM}$ - see details of the domain decomposition in
\cite{BL1}.

 Let the finite element geometry $\Omega$ is generated for any of  MM models
  of section \ref{sec:MM}
 and  let $\varepsilon = 1, \sigma =0$ in $\Omega_{\rm hyb} \setminus  \Omega$.
  Let  $\GammaT := \Gamma_{\rm hyb} \times (0,T)$  where  $\Gamma_{\rm hyb}$ is the boundary of the domain
 $\Omega_{\rm hyb}$.  In this case the model equation \eqref{forweq1} transforms to the wave equation in
  $\Omega_{\rm hyb} \setminus  \Omega$. For applications it is convenient to consider the first order absorbing boundary conditions \cite{EM} at 
  $\GammaT$ given by
\begin{equation}
 \pn E = - \pt E \text{  on }  \GammaT. \label{forwboundcond}
\end{equation}

\begin{figure}[tbp]
 \begin{center}
   \begin{tabular}{cc}
  \includegraphics[trim = 10.0cm 4.0cm 10.0cm 4.0cm, scale=0.2, clip=]{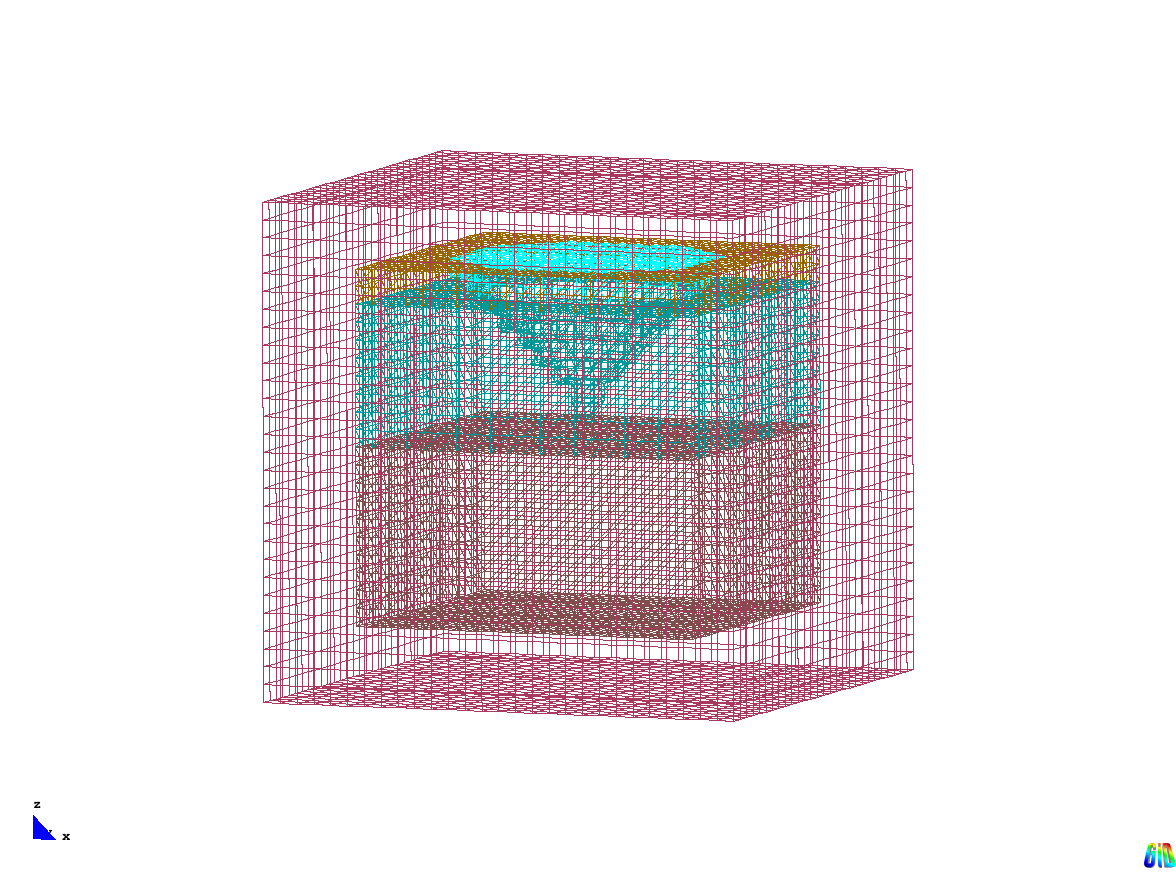} &
     \includegraphics[trim = 10.0cm 3.0cm 10.0cm 3.0cm, scale=0.2, clip=]{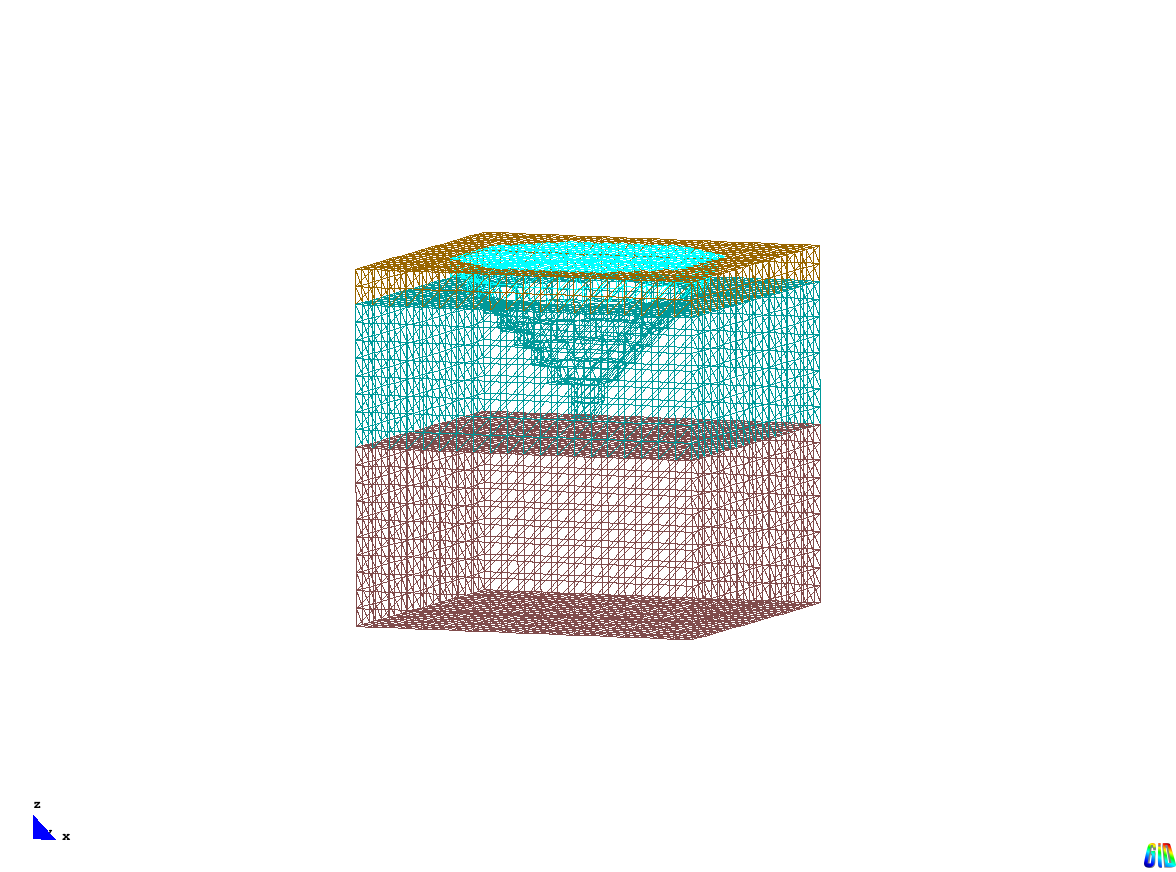}
\end{tabular}
\end{center}
 \caption{\textit{
 The domain  $\Omega_{\rm hyb} : = \Omega_{\rm FEM} \cup \Omega_{\rm FDM}$
in the DDFE/FDM used in computations. 
 Left: the combined finite element/finite difference  domain    $\Omega_{\rm hyb}$. Right: the finite element domain  $\Omega_{\rm FEM}$  corresponding to the model of MM at
 month 22
  with depth $d=4.06$ - see Figure \ref{fig:model3D} and Table  \ref{geotab}.}}
\label{fig:meshes}
\end{figure}

\section{Numerical study}

\label{sec:num}

\begin{figure*}
    \centering
    \begin{tabular}{c}
    \includegraphics[trim={3.4cm 8cm 1.5cm 6cm},clip, width=.9\textwidth]{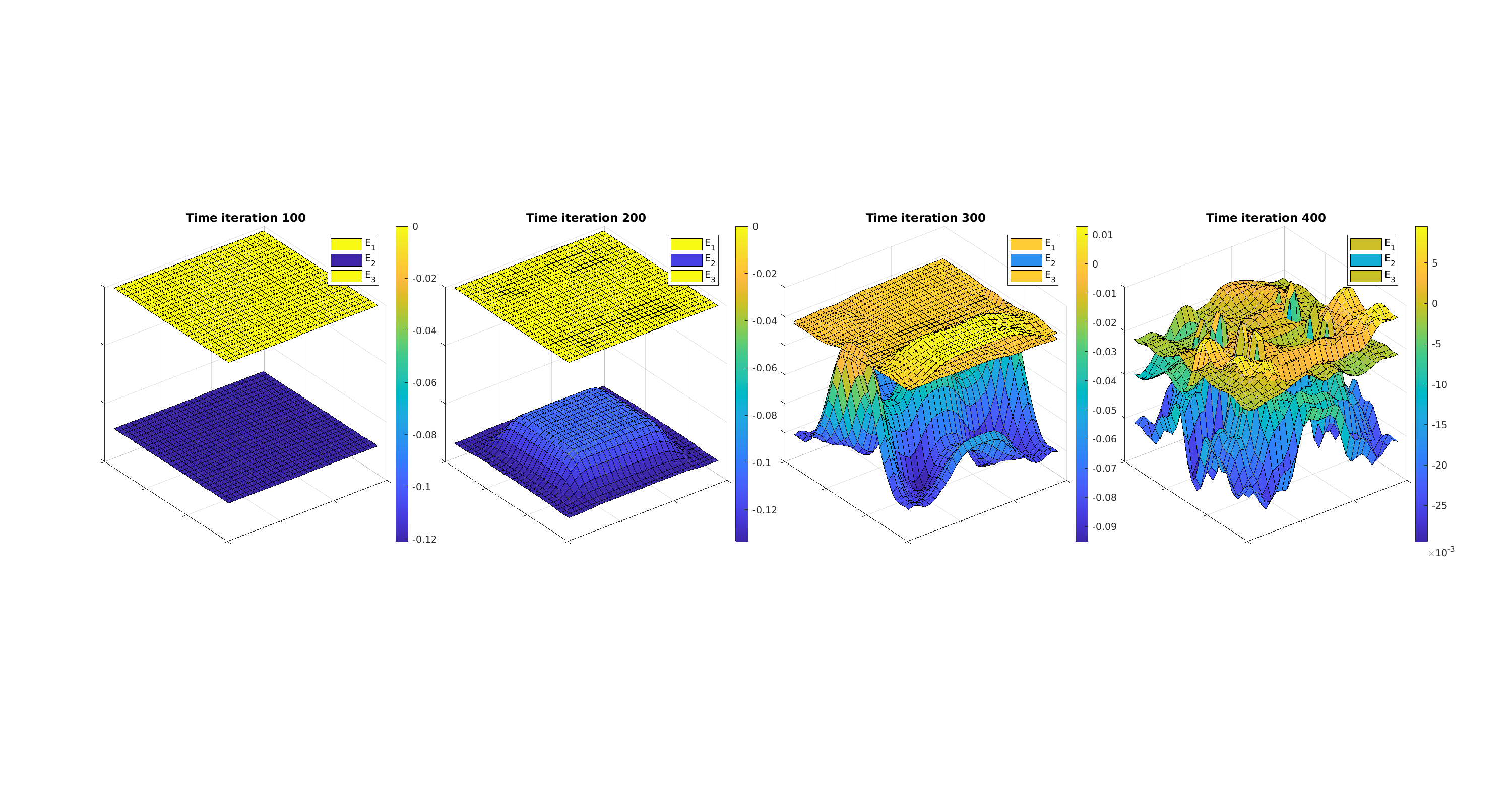} \\
    \begin{tabular}{c c c c}
    \includegraphics[trim={3cm 0cm 2cm 2cm},clip, width=.2\textwidth]{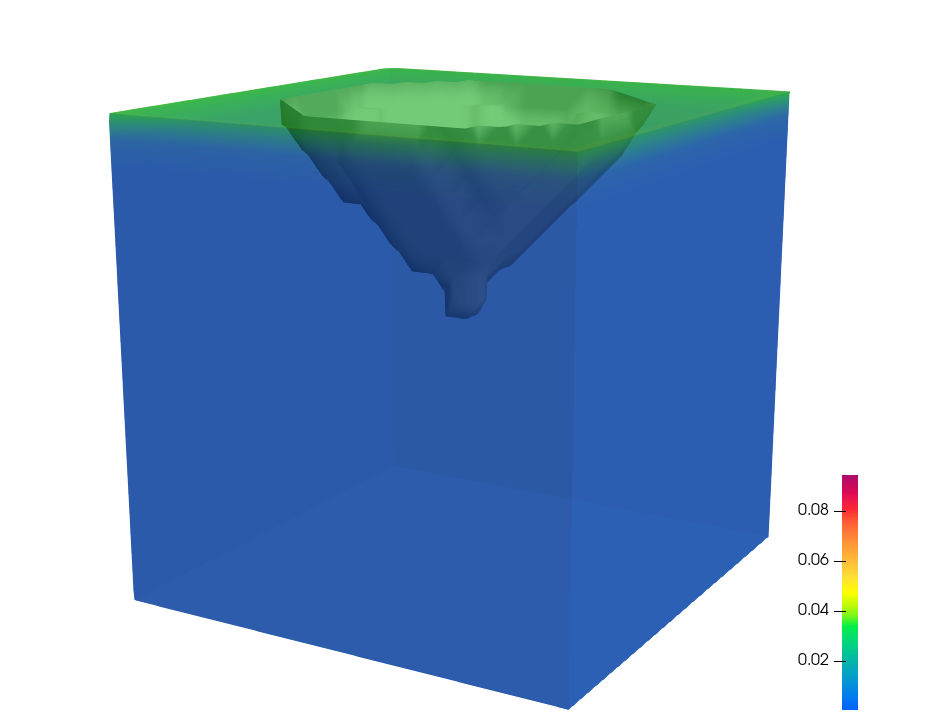} &
    \includegraphics[trim={3cm 0cm 2cm 2cm},clip, width=.2\textwidth]{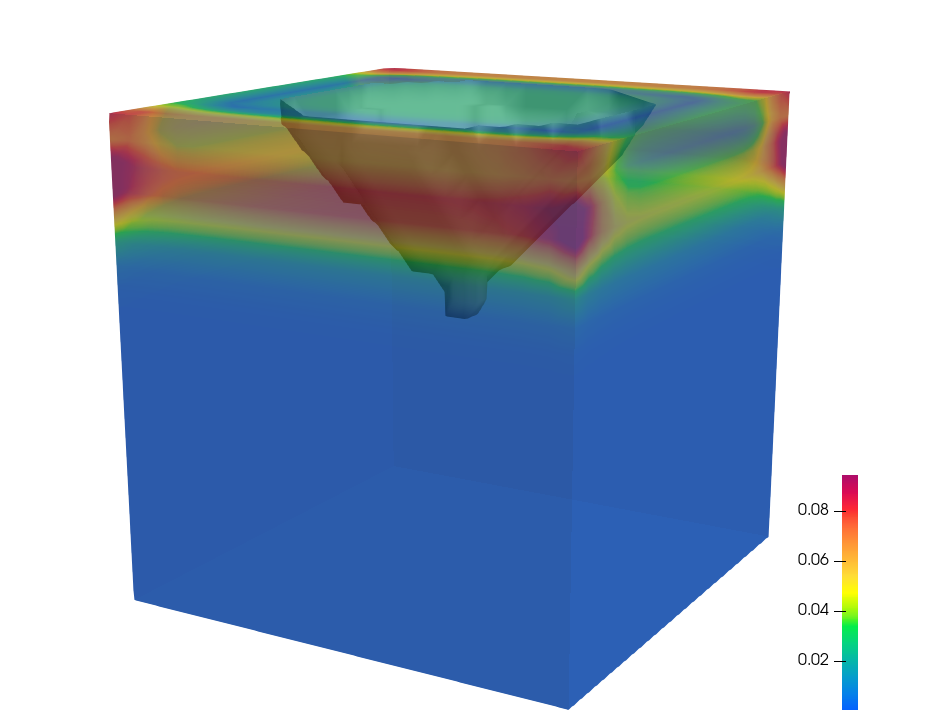} &
    \includegraphics[trim={3cm 0cm 2cm 2cm},clip, width=.2\textwidth]{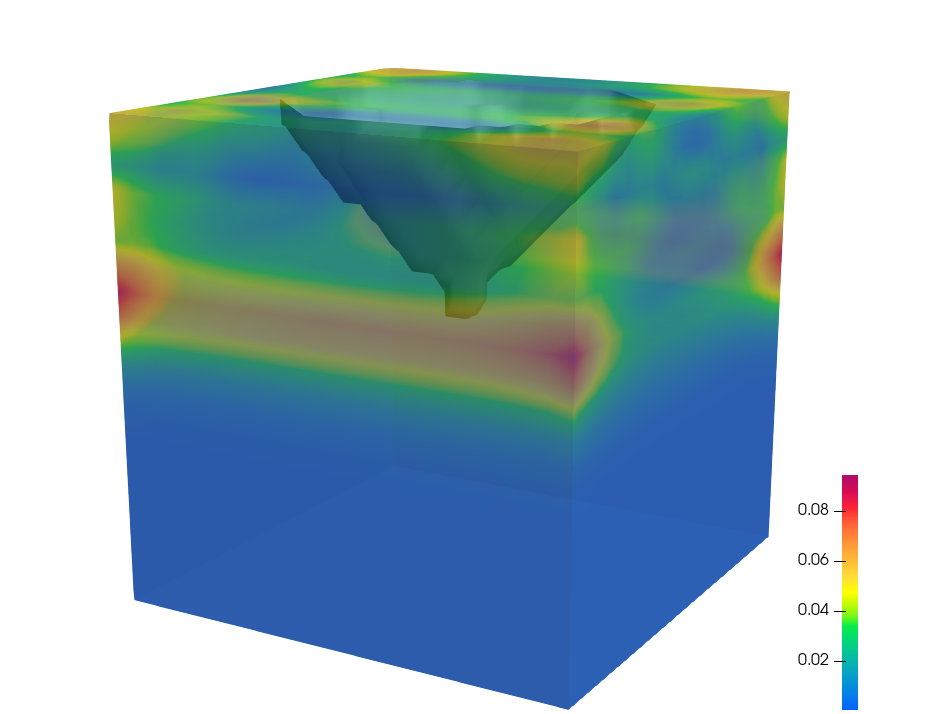} &
    \includegraphics[trim={3cm 0cm 2cm 2cm},clip, width=.2\textwidth]{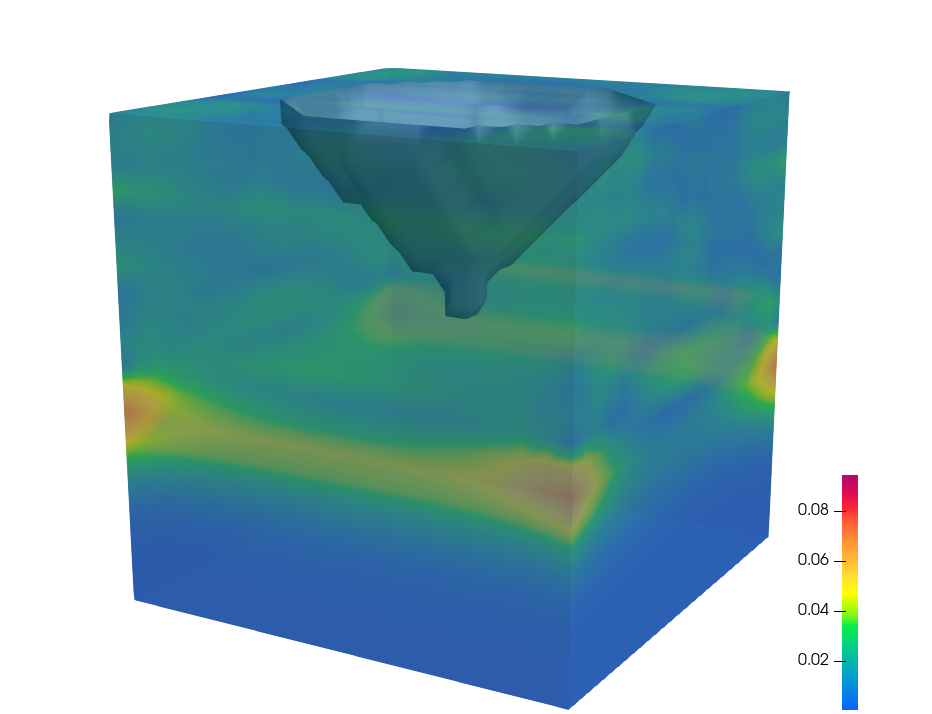}
    \end{tabular}
    \end{tabular}
    \caption{\textit{The top row shows the backscattered field measured at the top boundary of $E$. In the bottom row are images of the computed electric field $|E_h|$ which is shown as an opaque 3D rendering of the planar wave travelling through the domain $\Omega_{\rm FEM}$. Both rows correspond to MM at month 22 (with depth of tumor $d = 4.06$).}}
    \label{fig:forward}
\end{figure*}

This section demonstrates how the finite element 3D mesh
  for the MM model  with a depth  of $4.06$ mm
   (shown on Figure \ref{fig:meshes3D} and  Table  \ref{geotab})
is  used for generation of time-dependent
 backscattered data using the DDFE/FDM method of \cite{BL1} for the
 model equation \ref{forweq1} -- \ref{forwinitcond} with boundary
 condition \eqref{forwboundcond}.

\begin{figure*}[h!]
 \begin{center}
   \begin{tabular}{cc}
  \includegraphics[trim = 1cm 0.0cm 1cm 4cm, scale=0.2, clip=]{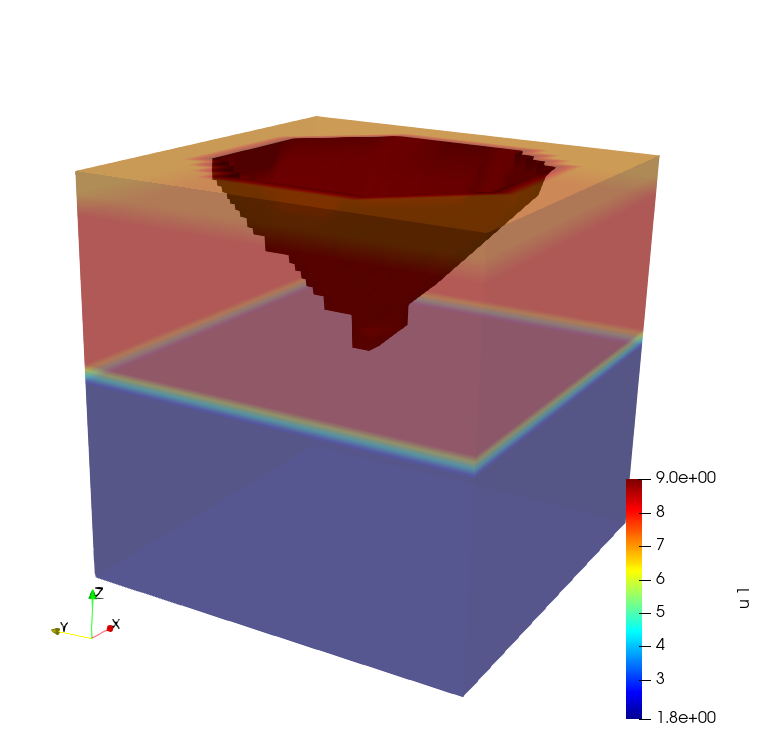 } &
     \includegraphics[trim = 1cm 0.0cm 1cm 4cm, scale=0.2, clip=]{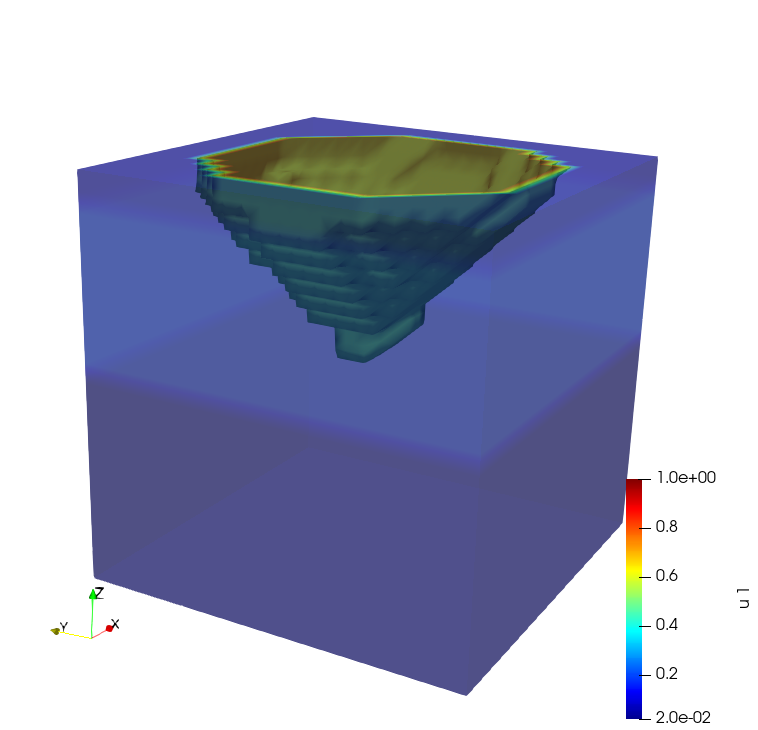 }\\
   $\max_\Omega \varepsilon_{r}/5 \approx 9$  &    $\max_\Omega \sigma/5 \approx 1.2$  \\
\end{tabular}
\end{center}
 \caption{\textit{Realistic weighted dielectric properties of MM  in skin at 6 GHz  corresponding to the MM at month 22 with depth $d=4.06$, see Tables  \ref{geotab},  \ref{tab:table1}, and  utilized for generation of time-dependent backscattered data via the DDFE/FDM method of \cite{BL1}. }}
\label{fig:numex1}
\end{figure*}

We choose the dimensionless computational domain
 $\Omega_{\rm hyb} : = \Omega_{\rm FEM} \cup \Omega_{\rm FDM}$
 (shown on Figure \ref{fig:meshes})
such that
$ \Omega_{\rm hyb}  = \left\{ x= (x_1,x_2, x_3) \in (-2, 12) \times (-2, 12) \times (-2,
 12) \right\}$,
and the computational domain $\Omega_{\rm FEM}$ as
$ \Omega_{\rm FEM} =
 \left\{ x= (x_1,x_2, x_3) \in (0, 10) \times (0, 10) \times (0,
 10) \right\}$.
The domain $\Omega_{\rm FEM}$ corresponds to the 3D model of MM
 with depth $d= 4.06$ mm, 
of the size $10 \times 10 \times 10$ mm,  shown on Figure \ref{fig:model3D}  for month 22.  Finite element discretization of this domain is presented on Figure \ref{fig:meshes}.
The values of $\varepsilon$
and $\sigma$ are assigned in $\Omega_{\rm FEM}$ accordingly to the weighted values of the Table
\ref{tab:table1}, and we set $\varepsilon =1$
and $\sigma = 0$ in $\Omega_{\rm hyb}  \setminus \Omega_{\rm FEM}$.
Figure \ref{fig:numex1} presents exact values of the relative dielectric permittivity and conductivity functions at 6 GHz in $\Omega_{\rm FEM}$. Further,  we use the same computational set-up as was used in \cite{ICEAA2025_KLB}.

\begin{table}[ht]
    \centering
    \begin{tabular}{| p{13em} | p{1cm} | p{1cm} | p{1cm} | }
    \hline
                &   &  &  depth   \\
    Tissue type  & $\varepsilon_r$  & $\sigma$  (S/m)  & (mm) \\
    \hline
    Immersion medium       & 32 & 4  &  2   \\
    \hline
    epidermis              & 35 & 4    & 1  \\
    \hline
     dermis                & 40 &  9  & 3.5  \\
    \hline
    Fat                    & 9  &  1   &  5.5  \\ 
    \hline
    Tumor stage 1          & 45 &  5   & $< 1$  \\
    \hline
    Tumor stage 2          & 50 &   5   & $ > 1$  \\
    \hline
    Tumor stage 3          & 60 &  6    &  $> 1$  \\
    \hline
    \end{tabular}
    \vspace{.1cm}\caption{\textit{Tissue types and corresponding  realistic values of $\varepsilon_r$  and $\sigma$  (S/m) at 6Ghz for skin with melanoma used in our numerical experiments.}}
    \label{tab:table1}
\end{table}

Figure \ref{fig:forward} demonstrates the computed scattered electric field $|E|$ of the model problem \ref{forweq1} -- \ref{forwboundcond} in the finite element domain at different times using DDFE/FDM method on the geometry
 shown on Figure \ref{fig:meshes}.

\begin{figure*}[tbp]
    \centering
    \includegraphics[trim={1cm 2cm 0cm 1cm},clip,width=0.95\linewidth]{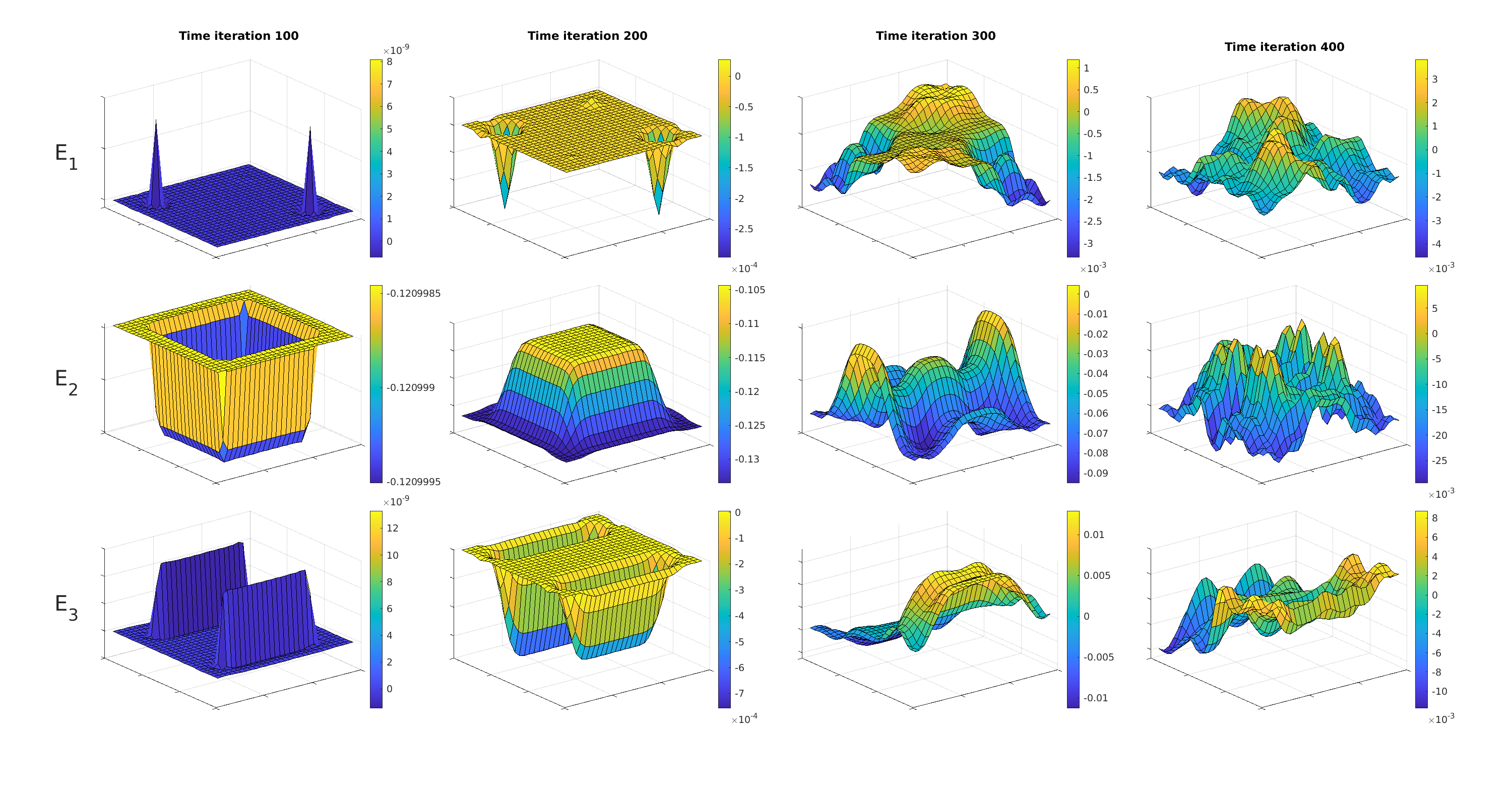}
    \caption{\textit{ The backscattered data $E_h = ({E_1}_h, {E_2}_h,{E_3}_h)$ measured on
    the backscattered boundary of the domain $\Omega_{\rm FEM}$  with
 MM at month 22 with depth d= 4.06
at different times.  }}
    \label{fig:backscatdata}
\end{figure*}

Figure \ref{fig:backscatdata} shows simulated backscattered data for
  the model problem \ref{forweq1} -- \ref{forwboundcond} at different
  times. This data can be used for testing of reconstruction
  algorithms to determine dielectric properties of the MM model with
  depth = 4.06 presented in Figure \ref{fig:numex1}.

\section{Conclusions}\label{S5}

\label{sec:concl}

In this study, we constructed finite element meshes for 3D models
simulating the growth of realistic malignant melanoma (MM) on the skin.
 These
 models are integrated with accurate representations of skin
 properties corresponding to the regions where MM develops.  To
 create the finite element meshes the software package WavES \cite{waves} was used, where
 we recently combined proposed geometrical models of
 melanoma growth from \cite{IEEE2024} with skin
 models investigated in \cite{ICEAA2024_BEN}.  Our numerical  tests
   show how  proposed 3D finite element meshes can be used to generate
 backscattered data for testing of
  reconstruction algorithms for determining the dielectric properties of the presented
   MM models.
All constructed finite element meshes
 as well as time-dependent data
are available upon request from the authors.

\section*{Acknowledgment}

The research of authors  is supported by the Swedish Research Council grant VR 2024-04459
 and STINT grant MG2023-9300.


\vspace{12pt}


\begin{thebibliography}{00}



\bibitem{ref1}  O. Abuzaghleh, B. D. Barkana, and M. Faezipour, “Noninvasive real-time
automated skin lesion analysis system for melanoma detection and pre-
vention,” IEEE J. Trans. Eng. Health Med., vol. 3, 2015, Art. no. 4300212,
doi: 10.1109/JTEHM.2015.2419612.
  
\bibitem{BakKok} A. B. Bakushinsky and M. Yu. Kokurin, \emph{Iterative
  Methods for Approximate Solution of Inverse Problems}, Springer,
  Dordrecht, The Netherlands, 2004.


\bibitem{BK} L. Beilina and M. V. Klibanov, \emph{Approximate global
  convergence and adaptivity for Coefficient Inverse Problems},
  Springer, New York, 2012.




\bibitem{BR1} L. Beilina, V. Ruas, On the Maxwell-wave equation coupling problem and its explicit finite element solution, \emph{Applications of Mathematics}, Springer, \url{https://doi.org/10.21136/AM.2022.0210-21}, 2022






 \bibitem{BL1}   L.~Beilina, E.~ Lindström, An Adaptive Finite Element/Finite Difference Domain Decomposition Method for Applications in Microwave Imaging, Electronics 2022, 11(9), 1359; https://doi.org/10.3390/electronics11091359

 \bibitem{BL2} L. Beilina,  E. Lindstr\"om,  A  posteriori error estimates and adaptive error control  for permittivity reconstruction in conductive media. In \emph{Gas Dynamics with Applications in Industry and Life Sciences}, Series: Springer Proceedings in Mathematics \& Statistics, Springer, PROMS, vol.429,  Cham (2023) 

\bibitem{ICEAA2024_BEN} L. Beilina, A. Eriksson and N. Neittaanmäki, "Frequency inversion method and device for malignant melanoma detection using RF/microwaves," 2024 International Conference on Electromagnetics in Advanced Applications (ICEAA), Lisbon, Portugal, 2024, pp. 794-799, doi: 10.1109/ICEAA61917.2024.10701729.

\bibitem{IEEE2024} J. Boparai, R. Tchinov, O. Miller, Y. Jallouli and M. Popović, "Models of Melanoma Growth for Assessment of Microwave-Based Diagnostic Tools," in IEEE Journal of Electromagnetics, RF and Microwaves in Medicine and Biology, vol. 8, no. 3, pp. 305-315, Sept. 2024, doi: 10.1109/JERM.2024.3430315.

 \bibitem{28} C. A. Bauman, P. Emary, T. Damen, and H. Dixon, “Melanoma in situ: A
case report from the patient’s perspective,” J. Can. Chiropractic Assoc.,
vol. 62, no. 1, pp. 56–61, 2018.

\bibitem{29} J. Beer, L. Xu, P. Tschandl, and H. Kittler, “Growth rate of melanoma
in vivo and correlation with dermatoscopic and dermatopathologic find-
ings,” Dermatol. Practical Conceptual, vol. 1, no. 1, pp. 56–67, 2011,
doi: 10.5826/dpc.0101a13.

\bibitem{30}  W. E. Damsky, L. E. Rosenbaum, and M. Bosenberg, “Decoding
melanoma metastasis,” Cancers (Basel), vol. 3, no. 1, pp. 126–163,
Dec. 2010, doi: 10.3390/cancers3010126.



\bibitem{EM}  B. Engquist  and A. Majda, Absorbing boundary conditions for
the numerical simulation of waves, \emph{\ Math. Comp.}, 31, 629-651, 1977.


\bibitem{31} R. Shayan, M. G. Achen, and S. A. Stacker, “Lymphatic vessels in
cancer metastasis: Bridging the gaps,” Carcinogenesis, vol. 27, no. 9,
pp. 1729–1738, Sep. 2006, doi: 10.1093/carcin/bgl031.

\bibitem{BuKr}
Maya de Buhan and Marie Kray,
A new approach to solve the inverse scattering problem for waves:
combining the {TRAC} and the adaptive inversion methods, Inverse Problems, 29(8), 2013
   



\bibitem{noninvasive} Fink C, Haenssle HA. Non-invasive tools for the diagnosis of cutaneous melanoma. \emph{Skin Res Technol.} 2017 Aug;23(3):261-271. doi: 10.1111/srt.12350.
  Epub 2016 Nov 22. PMID: 27878858.
  
 \bibitem{skindepth} Gershenwald JE, Scolyer RA, Hess KR et al.  Melanoma staging: Evidence-based changes in the American Joint Committee on Cancer eighth edition cancer staging manual. CA Cancer J Clin. 2017 Nov;67(6):472-492 

\bibitem{Ghavent} G. Chavent, \emph{Nonlinear Least Squares for Inverse Problems. Theoretical Foundations and Step-by-
Step Guide for Applications}, Springer, New York, 2009.

  \bibitem{GG} Gleichmann, Yannik G. and Grote, Marcus J.,
Adaptive Spectral Inversion for Inverse Medium Problems,
\emph{Inverse Problems},  39(12), 2023.  DOI: \url{10.1088/1361-6420/ad01d4}
  

\bibitem{gonch1} A. V.  Goncharsky, S. Y. Romanov,
  A method of solving the coefficient inverse problems of wave
tomography, \emph{Comput. Math. Appl.}, 2019;77:967–980.

\bibitem{gonch2} A. V. Goncharsky, S. Y. Romanov, S. Y. Seryozhnikov,
  Low-frequency ultrasonic tomography: math-
ematical methods and experimental results. Moscow University Phys Bullet. 2019;74(1): 43–51.

\bibitem{tomography} Editor(s):Pierre Grangeat, \emph{Tomography}, Wiley, 2009, DOI:10.1002/9780470611784


 



\bibitem{convex1} Vo Anh Khoa, Grant W. Bidney, Michael V. Klibanov,
  Loc H. Nguyen, Lam H. Nguyen, Anders J. Sullivan \& Vasily
  N. Astratov (2021), An inverse problem of a simultaneous
  reconstruction of the dielectric constant and conductivity from
  experimental backscattering data, \emph{Inverse Problems in Science and
  Engineering}, 29:5, 712-735, DOI: 10.1080/17415977.2020.1802447


\bibitem{TBKF2} N. T.  Th\'anh, L. Beilina, M. V. Klibanov,
  M. A. Fiddy, Imaging of Buried Objects from Experimental
  Backscattering Time-Dependent Measurements using a Globally
  Convergent Inverse Algorithm, \emph{SIAM Journal on Imaging
    Sciences}, 8(1), 757-786, 2015.

\bibitem{T}  A. N. Tikhonov,  A. V. Goncharsky,
 V. V. Stepanov and A. G. Yagola,
\emph{Numerical Methods for the Solution of Ill-Posed Problems}, London,
Kluwer, 1995.


\bibitem{itojin} K. Ito, B. Jin, \emph{Inverse Problems: Tikhonov theory and algorithms}, Series on Applied Mathematics, V.22, World Scientific,  2015.

  

\bibitem{ieee1} W.T. Joines, Y. Zhang, C. Li, and R. L. Jirtle, The
  measured electrical properties of normal and malignant human tissues
from 50 to 900 MHz', \emph{Med. Phys.}, 21 (4), pp.547-550, 1994.


\bibitem{KSS} S. Kabanikhin, A. Satybaev, and M. Shishlenin,
 \emph{Direct Methods of Solving Multidimensional Inverse Hyperbolic Problems}, VSP, Ultrecht, The Netherlands, 2004.

 


\bibitem{KBKSNF} AV Kuzhuget, L Beilina, MV Klibanov, A Sullivan, L Nguyen, MA Fiddy, Quantitative image recovery from measured blind backscattered data using a globally convergent inverse method,
IEEE transactions on geoscience and remote sensing 51 (5), 2937-2948, 2012

\bibitem{ICEAA2025_KLB}  G. Kyhn, E. Lindstr\"om, L. Beilina,  Reconstructing the dielectric properties of melanoma in 3D using real-life melanoma model, to appear in 
 2025 International Conference on Electromagnetics in Advanced Applications (ICEAA).

\bibitem{LB3} Lindström, E., Beilina, L. Energy norm error estimates and convergence analysis for a stabilized Maxwell’s equations in conductive media. Appl Math 69, 415–436 (2024).
\url{https://doi.org/10.21136/AM.2024.0248-23}

\bibitem{LB4} E. Lindström and L. Beilina, "A hybrid finite element/finite difference method for reconstruction of dielectric properties of conductive objects," 2024 International Conference on Electromagnetics in Advanced Applications (ICEAA), Lisbon, Portugal, 2024, pp. 788-793, doi: 10.1109/ICEAA61917.2024.10701914.

\bibitem{Pastorino} M. Pastorino, \emph{Microwave Imaging}, John Wiley \& Sons, Hoboken, NJ, 2010.

\bibitem{waves}  WavES, software package, \url{waves24.com}
\bibitem{wisconsin} E. Zastrow, S. K. Davis, M. Lazebnik, F. Kelcz,
  B. D. Veen, S. C. Hageness, Online repository of 3D Grid Based
  Numerical Phantoms for use in Computational Electromagnetics
  Simulations, \url{https://uwcem.ece.wisc.edu/MRIdatabase/}

\end{thebibliography}
\end{document}